\documentclass[11pt,a4paper,reqno]{amsart} 

\usepackage{amssymb,graphicx}

\usepackage{amssymb}

\newcommand{\koniec}{\begin{flushright}  $\Box $ \end{flushright}}
\newtheorem{theo1}{Theorem} 
\newtheorem{theo}{Theorem}[section] 
\newtheorem{prop}[theo]{Proposition}  
\newtheorem{lemma}[theo]{Lemma}
\newtheorem{defi}[theo]{Definition}

\def\theequation{\thesection.\arabic{equation}}

\newcounter{mnotecount}[section]

\renewcommand{\themnotecount}{\thesection.\arabic{mnotecount}}

\newcommand{\mnote}[1]
{\protect{\stepcounter{mnotecount}}$^{\mbox{\footnotesize
$
\bullet$\themnotecount}}$ \marginpar{
\raggedright\tiny\em
$\!\!\!\!\!\!\,\bullet$\themnotecount: #1} }

\newcommand{\hook}{{\setlength{\unitlength}{11pt}   
                   \begin{picture}(.833,.8)
                   \put(.15,.08){\line(1,0){.35}}
                   \put(.5,.08){\line(0,1){.5}}
                   \end{picture}}}
\newcommand{\CP}{\mathbb{CP}}
\newcommand{\C}{\mathbb{C}}
\newcommand{\HH}{\mathbb{H}}

\newcommand{\PP}{\mathbb{P}}
\newcommand{\sym}{\mbox{Sym}}
\newcommand{\PPP}{\mathbb{P}}
\newcommand{\RP}{\mathbb{RP}}
\def\by{\bf{y}}
\def\bp{\bf{p}}
\def\bq{\bf{q}}
\def\br{\bf{r}}
\def\lambd{\kappa}
\def\bs{\bf{s}}
\def\tf{\tilde{f}}

\newcommand{\R}{\mathbb{R}}

\def\p{\partial}
\def\OO{\mathcal{O}}

\def\ov{\overline}

\def\be{\begin{equation}}
\def\ee{\end{equation}}
\def\F{{\Psi}}
\def\w{{\wedge}}

\def\bea{\begin{eqnarray}}
\def\eea{\end{eqnarray}}

\newcommand{\spp}{\mathbb{S}}

\def\ov{\overline}

\begin{document} \date{16 May 2019}
\title{Conics, Twistors, and anti--self--dual tri--K\"ahler metrics}
\author{Maciej Dunajski}
\address{Department of Applied Mathematics and Theoretical Physics\\ 
University of Cambridge\\ Wilberforce Road, Cambridge CB3 0WA, UK.}
\email{m.dunajski@damtp.cam.ac.uk}

\author{Paul Tod}
\address{The Mathematical Institute\\
Oxford University\\
Woodstock Road, Oxford OX2 6GG\\ UK.
}
\email{tod@maths.ox.ac.uk}
\maketitle
\begin{abstract}
We describe the range of the Radon transform on the space $M$ of irreducible
conics in  $\CP^2$ in terms of natural differential operators associated to the
$SO(3)$--structure on $M=SL(3, \R)/SO(3)$ and its complexification. 
Following \cite{moraru} we show that for any function $F$ in 
this range, the zero locus of $F$ is a four--manifold admitting an 
anti--self--dual conformal structure which contains three different
scalar--flat K\"ahler metrics. The corresponding twistor space ${\mathcal Z}$ admits a holomorphic fibration over $\CP^2$. In the special case  where ${\mathcal Z}=\CP^3\setminus\CP^1$ the twistor lines
project down to a four--parameter family of conics which form  triangular Poncelet pairs
with a fixed base conic.
\end{abstract}
\section{Introduction}
The twistor construction of Penrose \cite{penrose}, and its Riemannian version
developed by Atiyah--Hitchin--Singer \cite{AHS} give one-to-one correspondences between anti--self--dual conformal structures $[\gamma]$ on a four--manifold $X$, and complex three--folds 
${\mathcal Z}$ with a four parameter family of rational curves. 
More conditions need to be imposed on ${\mathcal Z}$ if the conformal structure is to contain
a Ricci-flat metric. In this case there exists a holomorphic fibration
${\mathcal Z}\rightarrow \CP^1$ with a twisted symplectic form on the fibres.
A weaker condition is needed \cite{pontecorvo} if there exists a K\"ahler metric in $[\gamma]$. Then 
${\mathcal Z}$ admits an anti--canonical divisor given by a section of ${\kappa_{\mathcal Z}}^{-1/2}$,
where $\kappa_{\mathcal Z}$ is the holomorphic canonical bundle 
of ${\mathcal Z}$.
While there exist many explicit examples of twistor correspondences in both Ricci--flat \cite{H79, TW79}, and 
K\"ahler \cite{L91} cases, the resulting metrics in the conformal class are special in that they admit continuous groups of isometries. There are some notable exceptions in the Ricci--flat case, if
the twistor space fibers holomorphically over the total space of a line bundle 
${\mathcal O}(k)$. The corresponding hyper-K\"ahler metrics arise from the generalised Legendre
transform \cite{IR96}, and admit tri--holomorphic Killing spinors \cite{DM}, but in general no Killing vectors. 

Much less is known about the K\"ahler case. Moraru
\cite{moraru1, moraru} 
considered twistor spaces which holomorphically fiber over $\CP^2$.
He identified a set of second--order linear operators on the space $M$ of irreducible plane conics in $\CP^2$, and showed how any function $F$ in the kernel of these operators
gives rise to a conformal structure on the hypersurface $F=0$ in $M$. Moraru's papers do not contain
explicit examples, and an attempt to find such examples led us to this work. We show that some of 
Moraru's operators on $M$ are redundant, and we construct a set of independent operators
out of natural geometric structures on $M$: 
$M$ is an irreducible symmetric space carrying
an Einstein metric of  negative scalar curvature, and
a symmetric cubic three--form which, together with the  metric, gives an $SO(3)$ structure in the sense of \cite{bob}, \cite{frie}. The corresponding operators characterise the range of the Penrose--Radon transform
defined below on the space of conics.

We find several explicit functions in the range of this conic Penrose--Radon transform, and
construct corresponding conformal structures $[\gamma]$. There is a preferred metric, which we call
the barycenter metric, in each $[\gamma]$. The barycenter metric admits three linearly independent
solutions to the conformal Killing--Yano equations. Each such solution gives rise to an explicit conformal factor which reswcales the barycenter metric to a K\"ahler metric. K\"ahler metrics in a given ASD conformal class $(X, [\gamma])$ correspond \cite{DTkahler} to parallel sections of a certain connection on a rank ten vector bundle 
\[
E= {\Lambda^2}_+(X)\oplus\Lambda^1(X)\oplus {\Lambda^2}_-(X).
\]
The ASD conformal structures arising from our construction
admit a three--dimensional space of parallel sections of $E$. We call them 
tri--K\"ahler. 

In the next Section we shall introduce a $GL(2)$--structure isomorphism 
\[
\rho:T_m M\otimes \C\rightarrow\mbox{Sym}^4(\C^2)
\]
on the space
of irreducible real conics with no real points $M=SL(3, \R)/SO(3)$, and use it 
(Proposition \ref{prop_ein}) to construct an 
Einstein metric and a cubic three-form on $M$. They will  both be given in terms
of the bi--linear pairing $< , >_k: V_l\times V_n\rightarrow V_{l+n-2k}$ given by (\ref{k_th_trans}),
where $V_n=\mbox{Sym}^n(\C^2)$. In \S\ref{sec_radon}
these structures give rise
to  a metric Laplacian  $\Delta_g$, and another second order differential operator $\Box$ with values in $\Lambda^1(M)$ which characterise (Theorem \ref{theo_radon}) the range of the Penrose--Radon transform on conics as follows:
if $U$ is a neighbourhood of an irreducible conic  $C_m\subset \CP^2$ corresponding to $m\in M$, and  $f\in H^1(U, {\mathcal O}_{\CP^2}(-1))$.
Then the Penrose--Radon transform $F:M\rightarrow \R$ of $f$ satisfies the system of PDEs
\be
\label{theoremb_eq}
\Delta_g F=-\frac{1}{12}F, \quad \Box F=\frac{1}{24}d F.
\ee
Theorem \ref{theo_radon} makes some of the differential operators introduced in \cite{moraru} redundant, and clarifies the geometric meaning of the remaining operators.

In \S\ref{sec_asd}
we introduce a four--manifold $X$ as a hyper-surface in $M$ corresponding to the zero set of any function in the range of the Penrose--Radon transform.
We show (Theorem \ref{theo_1_forms} and Theorem \ref{theo_2_forms}) that this conformal structure
is anti--self--dual, and that it contains a barycentre metric.
We shall give several explicit examples
of barycentre metrics, including a Ricci--flat example.
Finally 
in \S\ref{quad_map_s}
we consider a holomorphic fibration of the complement of a rational normal curve
in $\CP^3$ over $\CP^2$. We characterise (Proposition \ref{theo_ponce}
and Proposition \ref{propFG}) the four--parameter family of conics which arise
as images of lines in $\CP^2$ under a quadratic map. This brings up some 
 classical 19th century  projective geometry involving loci of Gergonne points
and the Poncelet porism \cite{Gallaty, cayley} . We shall prove
\begin{theo1}
\label{theo3_intro}
Let $Q:\CP^3\setminus\CP^1\rightarrow \CP^2$ be the projectivisation of
the moment map for the symplectic
$SL(2, \C)$ action on the space of cubics with at least two distinct roots.
The image of lines in $\CP^3\setminus\CP^1$ are conics in $\CP^2$ corresponding to a hypersurface
$X\subset M$. If $x\in X$, then there exists a triangle inscribed in 
the corresponding conic ${C_x}$ and circumscribed about
a base conic $(Z^1)^2+(Z^2)^2+(Z^3)^2=0$. Equivalently
\[
2\mbox{Tr}({A_x}^2)-\mbox{Tr}(A_x)^2=0,
\]
where ${A_x}={A_x}^T$ is the non--singular symmetric 3 by 3 matrix defining the conic ${C_x}$.
\end{theo1}
\subsection*{Acknowledgements.} We are grateful 
to Claude LeBrun and Lionel Mason for discussion about twistor theory,
and to Robert Bryant, Nigel Hitchin and Miles Reid
for correspondence clarifying the quadratic map discussed in 
\S\ref{quad_map_s}. The work of M.D. has been partially supported by STFC consolidated grant no. ST/P000681/1. Part of this work was done while P.T. held the 
Brenda Ryman Visiting Fellowship in the Sciences at Girton College, Cambridge, and he gratefully acknowledges the hospitality of the College.

\section{The space of plane conics}
\subsection{$GL(2)$ structure and twistor theory}
Let  $Z=[Z^1, Z^2, Z^3]$
be the homogeneous coordinates in $\CP^2$.
A general conic in $\CP^2$ is of the form
\be
\label{homogeneous}
f([Z], A)\equiv ZAZ^T=0,
\ee
where $A$ is a complex symmetric matrix.
Thus the space of conics is  $\PP(\mbox{Sym}^2(\C^3))=\CP^5$. 
We shall consider the five--dimensional space $M_\C$  of 
irreducible conics normalised by $\mbox{det}(A)=1$.
The projective group $PSL(3, \C)$ acts on $M_\C$ transitively by
\be
\label{sl_action}
A\longrightarrow N A N^T
\quad\mbox{where}\quad
N\in SL(3, \C),
\ee
so $M_\C\cong SL(3, \C)/SO(3, \C)$ as $SO(3, \C)$ stabilises the conic
$A={\bf 1}$. 
The conic $A={\bf 1}$,  or equivalently $ZZ^T=0$ admits no real  points, so it belongs
to $\CP^2\setminus\RP^2$. Its $SL(3, \R)$ orbit consists of conics with real coefficients but  no real points.
The real five--dimensional  manifold of such conics is $M=SL(3, \R)/SO(3, \R)$. It is a real  slice in 
the complex manifold $M_\C$.

 The vector fields in $M_\C$ can be canonically identified
with homogeneous fourth order polynomials in two variables - this will play a 
role in what follows. To formalise it, let us first make a definition
\begin{defi}
A $GL(2)$ structure on an $(n + 1)$- complex dimensional complex manifold $M_\C$ 
is an isomorphism:
\be
\label{gl2def}
T M_\C \cong \mbox{Sym}^n(\spp)
\ee
where $\spp$ is a rank--two complex symplectic vector bundle over $M_\C$.
\end{defi}
\noindent
This structure was called a {\em paraconformal structure} in \cite{DT},
and $G_n$ structure in \cite{bryant}. In practice a $GL(2)$ structure is specified by a one--form
$S$ on $M_\C$ with values in $\mbox{Sym}^n(\C^2)$. Let $e^1, e^2, \dots, e^{n+1}$ be
$(n+1)$ independent one--forms on $M_\C$. Then
\[
S=t^n \;e^1+n t^{n-1}s\; e^2 +\frac{1}{2} n(n-1) t^{n-2}s^2\; e^3+\dots +s^n\; e^{n+1},
\]
and the isomorphism (\ref{gl2def}) is given by $\rho(V)=V\hook S$, where
$V\in \Gamma(TM_\C)$ and $\hook$ denotes contraction of a vector field with a form.

For even $n=2l$ there are two real forms of the $GL(2)$ structure. 
In the {\em indefinite} case the one forms $e^1, \dots, e^{n+1}$ can be taken to be real,
and $[s, t]\in \RP^1$. In the {\em positive definite} case
\[
\overline{S(\lambda)}=(-1)^l\ov{\lambda}^{2l}S(-1/\ov{\lambda}), \quad\mbox{where}\quad
\lambda=t/s.
\] 
In what follows we shall be interested in the case $n=4$, where the holomorphic five--manifold
$M_\C$ admits real form $M$ with a positive definite quadratic form. In this case
\be
\label{para_con_exp}
S[s, t]=t^4\; e^1+4 t^3s\; e^2 +6 t^2 s^2\; e^3-4 ts^3\; \ov{e^2}+
s^4 \ov{e^1},
\ee
where the one--forms $e^1, e^2$ are complex,  and $e^3$ is real.
\begin{prop}
\label{GL2prop}
The space of plane conics carries a $GL(2)$ structure.
\end{prop}
\noindent
{\bf Proof}.
All plane conics are rational curves and two neighbouring conics intersect at four points. Thus if 
$C_m\cong \CP^1$ is a conic corresponding to $m\in M_\C$
then the normal bundle $N(C_m)\cong \OO(4)$.
The obstruction group $H^1(\CP^1, \OO(4))=0$,
so $M_\C$ is a locally complete family, and there exists a
Kodaira isomorphism \cite{kodaira}
\be
\label{paracon}
T_m M_\C\cong H^0(C_m, N(C_m))=\mbox{Sym}^4(\C^2)
\ee
identifying vectors tangent to $M_\C$ with binary quartics. Therefore $M_\C$ carries a $GL(2)$ 
structure. 
The isomorphism (\ref{paracon}) is specified by a fourth-order polynomial 
$S$ homogeneous in two variables and
with values in $\Lambda^1(M_\C)$. Given such $S$,  the  binary quartic  
corresponding under (\ref{paracon}) to a vector field $V\in TM_\C$ is 
the contraction $V\hook S$.
\koniec
\subsubsection{Transvectants and Invariants}
\label{sub_transv} 
In the reminder of this section we shall use the $GL(2)$ structure on $M_\C$
to construct a conformal structure, and a symmetric cubic three--form. Both 
structures originate from classical invariants of binary quartics 
\cite{Elliot}, and these invariants can be conveniently introduced using 
the transvectant/spinor notation together with some representation theory of
$\mathfrak{sl}(2, \C)$.

Let $V_{l}=\mbox{Sym}^{l}(\C^{2})$ be the $(l+1)$--dimensional complex vector space of binary quantics of degree $l$. These quantics are the same as polynomials
in $[s, t]$ homogeneous of degree $l$.
Given two binary quantics $\phi\in V_n$ and $\psi\in V_l$, the $k$th transvectant is a map $<\;,\;>_k:V_l\times V_n\rightarrow
V_{l+n-2k}$ given by a quantic of degree $l+n-2k$
\be
\label{k_th_trans}
<\phi, \psi>_k:=\sum_{j=1}^k(-1)^j{k\choose j}\frac{\p^k \phi}{\p s^{k-j}\p t^j}
\frac{\p^k \psi}{\p s^{j}\p t^{k-j}}.
\ee
This map is $SL(2, \C)$ equivariant, and
is symmetric for $k$ even and skew-symmetric for $k$ odd.
Where possible,  we will use the notation (\ref{k_th_trans}) instead of the 
spinor notation \cite{md_rp} which we summarise in the Appendix.  For example
if  $\phi\in V_3$ and $\psi\in V_4$ then $<\phi, \psi>_2$ is a cubic 
proportional to ${\phi_{(A}}^{DE}\psi_{BC)DE}$.

\subsection{Conformal structure and Einstein metric}
\label{sub22}
In the next Proposition we shall give a twistor construction of 
an Einstein metric on a five dimensional manifold  $M$ which is a real slice of $M_\C$. 
The metric itself is well known in the theory of symmetric spaces \cite{sym_spaces}. Our treatment will follow the twistor procedure of
\cite{DG}.
\begin{prop}
\label{prop_ein}
Let $[g_\C]$ be the conformal structure on the space of complex irreducible 
conics $M_\C=SL(3, \C)/SO(3, \C)$ such that 
the  vector field $V\in \Gamma(TM_\C)$ is null iff the corresponding 4th order 
polynomial has a vanishing transvectant $<V, V>_4=0$. 
\begin{itemize}
\item There exists
a real form $M=SL(3, \R)/SO(3)$ of $M_\C$ where $[g_\C]$ gives a positive definite conformal structure $[g]$, and such that there exists a 
Riemannain metric $g\in [g]$
which is Einstein, has negative scalar curvature, and admits
$SL(3, \R)$ as its isometry group. 
\item Local coordinates $(a, b, p, q, r)$ can be used on $M$ 
to parametrise conics, such 
that then
\be
\label{final_metric}
g=8[(2da+db)^2+ 3db^2+
e^{2a+4b}dp^2+e^{2a-2b}dr^2+e^{4a+2b}(dq-pdr)^2].
\ee
\end{itemize}
\end{prop}
\noindent
{\bf Proof.}
 A vector in $V\in \C\otimes TM$ is null
if the corresponding quartic has vanishing transvectant
i. e. 
\be
\label{con_con}
<V, V>_4=4\alpha\epsilon-8\beta\delta+6\gamma^2,\quad \mbox{where} \quad 
V=\alpha t^4+4\beta t^3s+6\gamma
t^2s^2+4\delta ts^3+\epsilon s^4.
\ee
This nullity condition can be interpreted in terms of
the roots of $V$: they can be transformed to vertices
of a regular tetrahedron in $\CP^1$ (if they are 
distinct), or at least three of the roots 
coincide\footnote{Such a conformal structure exists on any odd--dimensional
manifold with a $GL(2)$--structure. The geometric interpretation in terms
of roots of the associated homogeneous polynomial also exists, but is a little more complicated \cite{md_rp}. For even $n$, say $n=2l$, the real forms are as follows: in the  {\em indefinite}
case the quadratic form $<\rho(V), \rho(V)>_{2l}$ has signature $(l, l+1)$, and is invariant
under $SL(2, \R)$; in the {\em positve definite} case this quadratic form is positive--definite, and is invariant
under $SU(2)$.
} \cite{Elliot}.

To find explicit forms of the metric and the $GL(2)$ stucture on $M$ we use the fact that
any symmetric matrix $A$ is of the form 
\be
\label{upper}
A=BB^T,
\ee
where $B$ is upper
triangular, so that the general conic (\ref{homogeneous}) 
is $ZBB^T Z^T=0$. A convenient parametrisation of $B$ turns out to be
\[
B=\left(
\begin{array}{ccc}
e^c & pe^b& qe^a\\
0 & e^b &re^a\\ 
0 & 0 &e^a
\end{array}
\right ),
\]
where $(a, b, p, q, r)$ are real  coordinates on $M$.
To single out the irreducible conics we chose the normalisation 
$\mbox{det}(B)=1$ which corresponds to $a+b+c=0$. 

To find the paraconformal structure consider a neighbouring conic  where $B$ is replaced
by $B+dB$ in (\ref{homogeneous}). It will intersect
the conic (\ref{homogeneous}) at four points corresponding to the roots of a quartic
\be
\label{para_S}
S(s, t)=\sum_{\alpha,\beta} \frac{\p f}{\p B_{\alpha\beta}}|_{Z=Z([s, t], B)}dB_{\alpha\beta},
\ee
where $Z([s, t], B)$ is a rational parametrisation of the general conic and $\mbox{det}(B)=1$.
To find this rational parametrisation
consider a conic $W=[W^1, W^2, W^3]$ given by
\[
(W^1)^2+(W^2)^2+(W^3)^3=0
\] 
which is
(\ref{homogeneous}) with $B={\bf 1}$.
Any other conic is projectively equivalent to this one,
so we can use
(\ref{para_S})  with $f=ZBB^TZ^T$, and 
the parametrisation 
\begin{eqnarray}
\label{parametrisation}
Z&=&W B^{-1}, \quad\mbox{where}\quad  W=[s^2-t^2, 2ts, i(s^2+t^2)]\nonumber\\
&=&[e^{a+b}(s^2-t^2), 2e^{-b}st-e^{a+b}p(s^2-t^2),  e^{a+b}(pr-q)(s^2-t^2)-2re^{-b}st
+ie^{-a}(s^2+t^2)].
\nonumber
\end{eqnarray}
The formula (\ref{para_S}) leads to 
\[
S=W\Omega W^T, \quad\mbox{where}\quad \Omega\equiv B^{-1}dB+ 
(B^{-1}dB)^{T}.
\]
Note that the $SL(3, \R)$ action on $M$ takes the form $B\rightarrow N B$. This does not preserve
the upper--triangular matrices $B$, but it preserves $B^{-1}dB$. Therefore both the paraconformal structure encoded in the quartic $S$ and the resulting metric 
(\ref{con_con}) are $SL(3, \R)$ invariant. Computing $S$ explicitly  gives
\begin{eqnarray*}
S(s, t)&=&t^4(\Omega_{11}-\Omega_{33}-2i\Omega_{13})
+4t^3s(i\Omega_{23}-\Omega_{12})+6t^2s^2
(\frac{2}{3}\Omega_{22}-\frac{1}{3}\Omega_{11}-
\frac{1}{3}\Omega_{33})\\
&+&4ts^3(i\Omega_{23}+\Omega_{12})+s^4(\Omega_{11}-\Omega_{33}+2i\Omega_{13}).
\end{eqnarray*}
This leads to the positive definite  
$SL(3)$--invariant conformal structure (\ref{con_con})
\be
\label{metric_frame}
g=4(\Omega_{12}^2+\Omega_{13}^2+\Omega_{23}^2)+(\Omega_{11}-\Omega_{33})^2
+3(\Omega_{11}+\Omega_{33})^2
\ee
where we have used $\Omega_{11}+\Omega_{22}+\Omega_{33}=0$.
To find the explict form  of the conformal structure compute
\[
B^{-1}dB=\left(
\begin{array}{ccc}
-da-db & e^{2b+a}dp& e^{2a+b}(dq-pdr)\\
0 & db &e^{a-b}dr\\ 
0 & 0 &da
\end{array}
\right ).
\]
Using $(a, b, p, q, r)$ as real coordinates on $M$ yields the Riemannian 
metric
(\ref{final_metric}).
We verify by explicit calculation that this metric is Einstein 
with  negative scalar curvature equal to $-15/16$.
The isometry group of 
(\ref{final_metric}) is $SL(3, \R)$. The Killing vectors
generating this goup are given in the Appendix by the formula
(\ref{5d_killing}).
\koniec
\subsection{$SO(3)$ structure on the space of conics}
\label{subso3}
In this Section we shall reveal an additional structure on the space of
irreducible conics which will (in \S\ref{sec_radon}) play a role
in a characterisation of the Penrose--Radon transform on coincs.
\begin{defi}
An integrable $SO(3)$ structure on a five--dimensional manifold $M$
is a pair $(g, G)$, where $g$ is a Riemannain metric, and $G$ is a symmetric
three--form on $TM$ such that \cite{frie, bob}.
\be
\label{so3structure}
\nabla_a G_{bcd}=0, \quad
6{G^a}_{(bc}G_{de)a}=g_{(bc}g_{de)}, \quad G_{abc}g^{ab}=0.
\ee
\end{defi}
We have seen that the fourth transvectant endows the space of conics with an 
Einstein metric. We shall now see that another transvectant  operation
gives an $SO(3)$ structure on $(M, g)$.
\begin{prop}
\label{so3prop}
Let $g$ be the Einstein metric on the space of conics from Proposition \ref{prop_ein}, such that
$g(V, V)=<V, V>_4$. The symmetric cubic form
\[
G:TM\odot TM\odot TM\rightarrow \R
\]
given by $G(V, V, V)=<<V, V>_2, V>_4$ gives an integrable $SO(3)$ structure.
\end{prop}
\noindent
{\bf Proof.}
We will establish the identities (\ref{so3structure}) directly.
The $GL(2)$ structure $S\in 
\Lambda^1(M)\otimes\mbox{Sym}^4(\C^2)$ is given by
(\ref{para_con_exp}),
where $e^1, \dots, e^5$ is the pentad of one--forms on $M$ given by
\begin{eqnarray*}
e^1&=&-2(db+2da+ie^{2a+b}(dq-pdr)), \quad
e^2=-(e^{a+2b}dp-i e^{a-b}dr), \quad 
e^3=2db, \\
 e^5&=&\overline{e^1}, \quad
e^4=-\overline{e^2}.
\end{eqnarray*}
This  gives rise to the metric $g$ from Proposition \ref{prop_ein},  
and a symmetric three--form $G$ as
\begin{eqnarray}
\label{so3vector}
g&=& 2e^1 \odot e^5-8e^2\odot e^4+6 e^3\odot e^3\\
G&=&6(e^1\odot e^5\odot e^3 +2 e^2\odot e^4\odot e^3 -e^3\odot e^3\odot
e^3- e^1\odot e^4\odot e^4-e^5\odot e^2\odot e^2),\nonumber
\end{eqnarray}
where
\begin{eqnarray*}
u\odot v&:=&\frac{1}{2} \Big(u\otimes v+ v\otimes u\Big)\\
u\odot v \odot w&:=&\frac{1}{6}\Big(u\otimes v\otimes w+ u\otimes w\otimes v+
v\otimes u\otimes w +  v\otimes w\otimes u+
w\otimes u\otimes v+ + w\otimes v\otimes w\Big).
\end{eqnarray*}
The metric and the cubic form satisfy (\ref{so3structure}).
Moreover the $SL(3)$ action 
(\ref{sl_action})
on $M$ preserves $G$.
\koniec
\section{The Penrose--Radon transform on the space of conics}
\label{sec_radon}
Let ${\mathcal L}\rightarrow \CP^2$ be a holomorphic line bundle
such that ${\mathcal L}|_{C_m}\cong \OO_{\CP^1}(-2)$ for any conic and let
$f\in H^1(U, {\mathcal L})$, where $U$ is a neighbourhood of $C_m$ in $\CP^2$. Restricting $f$ to a conic $C_m$ and integrating over a contour
$\Gamma\subset C_m$ defines a function on $M$ by
\be
\label{contour_int}
F(m)=\oint_{\Gamma\subset C_m} f[Z^1(s, t, m), Z^2(s, t, m), Z^3(s, t, m)] (sdt-tds).
\ee
If $[Z^1, Z^2, Z^3]$ are three quadrics in $[s, t]$, then
\[
\theta=f(Z^{i})(sdt-tds), \quad d\theta=\Big(Z^{i}
\frac{\p f}{\p Z^{i}}+f\Big)ds\wedge dt.
\]
Thus  we require 
$f$ to be homogeneous of degree $-1$ in $Z$.
For example taking
\[
f=\frac{Z^2}{Z^1Z^3}
\]
and using a contour enclosing the singularity
at $[s, t]=[1, 0]$ in $\CP^1$ gives (on taking the real part)\footnote{Some other solutions are obtained as the residues:
\[
Z_3/{(Z_1)}^2\rightarrow \frac{q-pr}{e^{a+b}}, \quad
Z_2/(Z_1Z_3)\rightarrow \frac{re^a+ie^b}{e^b(r^2e^{2a}+e^{2b})}, \]
\[
Z_3/(Z_1Z_2)\rightarrow\frac{r-ie^{b-a}}{e^{a+b}}, \quad
Z_2/{(Z_1)}^2\rightarrow pe^{-a-b}, \quad 1/Z_1\rightarrow -e^{-a-b}.
\]
}
\be
F=\frac{e^{a-b}r}{e^{2b}+r^2 e^{2a}}.
\ee
\subsection{Differential operators on $SL(3, \R)/SO(3)$, and the 
range of the Penrose--Radon transfom}
We aim to characterise the range of the transform (\ref{contour_int}) as the kernel of 
some differential operators on $M$. One natural operator is the 
metric Laplacian 
$\Delta_g=g^{ab}\nabla_a\nabla_b$
of the metric (\ref{final_metric})
 explicitly given in the Appendix 
by (\ref{aplaplacian}).
Another operator can be defined for any $SO(3)$ structure as follows.
Define a second--order operator $\Box$ with values
in $\Lambda^1(M)$ by
\[
\Box F=({G_{a}}^{bc}\nabla_b\nabla_c F)e^a, 
\]
where the indices are raised and lowered by $g$ and its inverse.
\begin{theo}
\label{theo_radon}
Let the function $F:M\rightarrow \R$ belong to the image of the Penrose--Radon transform (\ref{contour_int}). Then 
\be
\label{moraru_eq}
\Delta_g F=-\frac{1}{12}F, \quad \Box F=\frac{1}{24}dF.
\ee
\end{theo}
\noindent
{\bf Proof.}
This must follow from the abstract machinery of the Penrose transform 
\cite{beastwood}, but here we present a concrete proof based on an explicit computation which hinges on 
a link between plane conics and the geometry of a certain fifth-order ODE.
We start with the integral  (\ref{contour_int})
 and transform it by observing that each $Z^i$ is homogeneous of degree 2 in the pair $[s,t]$, while $f$ is homogeneous of degree $-1$ in the $Z^i$. Thus, with $\lambda=t/s$,
\[F(m)=\oint f(Z^i(\lambda))d\lambda.\]
Introduce $x(\lambda)=Z^1(\lambda)/Z^3(\lambda),y(\lambda)=Z^2(\lambda)/Z^3(\lambda)$ to obtain
\[F(m)=\oint \tf(x,y)\frac{d\lambda}{Z^3(\lambda)},\]
where $\tf(x,y)=f(x,y,1)$.

The pair $(x(\lambda),y(\lambda))$ defines a conic in parametric form but the conic can also be defined by giving $y$ as a function of $x$. In this case, with overdot for $d/d\lambda$, we may calculate
\[p:=\frac{dy}{dx}=\frac{Z^3\dot{Z}^2-Z^2\dot{Z}^3}{Z^3\dot{Z}^1-Z^1\dot{Z}^3},\]
and then
\[q:=\frac{d^2y}{dx^2}=\frac{\Delta(Z^3)^3}{(Z^3\dot{Z}^1-Z^1\dot{Z}^3)^3},\]
where
\[\Delta=\ddot{Z^1}(Z^2\dot{Z^3}-Z^3\dot{Z^2})+\ddot{Z^2}(Z^3\dot{Z^1}-Z^1\dot{Z^3})+\ddot{Z^3}(Z^1\dot{Z^2}-Z^2\dot{Z^1}),\]
which can be calculated from \S\ref{sub22}, with the result $\Delta=8i$, in particular $\Delta$ is constant. It is  now easy to see that
\[\frac{d\lambda}{Z^3(\lambda)}=\frac{dx}{\dot{x}Z^3}=q^{1/3}\Delta^{-1/3}dx,\]
and so, up to a multiplicative constant, which we ignore
\be\label{eq23}F(m)=\oint \tilde{f}(x,y(x))q^{1/3}dx.\ee
This will be the form of the integral with which we calculate. Both $y(x)$ and $q$ depend on $m$, the conic on which the integral is performed. 
It is straightforward to see that (\ref{eq23}) agrees with the formula on p47 of \cite{moraru1}.

For the next stage we review and revise some theory from \cite{DT}. Consider the following differential equation for dependent variable $y(x)$:
\be\label{eq1}y^{(5)}=\Lambda(x, y, p, q, r, s)=-\frac{40}{9}\frac{r^3}{q^2}+5\frac{rs}{q}\ee
with $p=y',q=y'',r=y''',s=y''''$. It is well-known, but in any case easy to see, that the solutions of (\ref{eq1}) are the conics (first observe that the equation can be written as $(q^{-2/3})'''=0$).
Write the solution as
\[y=Z(x,X^{\bf a}),\quad{\bf a}=1,\ldots, 5,\]
so that $X^{\bf a}$ are coordinates on the solution space, that is on the space of conics, and concrete indices are bold. It will be convenient to choose the $X^{\bf{a}}$ to be the derivatives 
of $y$ from order zero to order four at some fixed but arbitrary choice of $x$. Call these $(\by,\bp,\bq,\br,\bs)$.

Following \cite{DT}, we suppose that there is a metric and spin-structure on the solution space such that the gradient of $Z$ is a quartic with a quadrupole root
\be
\label{special_spinor}
{\by}_{,a}:=Z_{,a}=\iota_A\iota_B\iota_C\iota_D,\ee
for some spinor $\iota_A$ (and abstract indices are italic). There is an abstract justification for this assumption, as the space of conics admits
the $GL(2)$ structure of Proposition \ref{GL2prop}. With prime 
for $d/dx$ we may suppose
\[\iota'_A=Po_A\]
for some $P$ to be found, where $o_A$ completes a spinor basis, so that $o_A\iota^A=1$ and $\epsilon_{AB}=o_A\iota_B-o_B\iota_A$. Now necessarily 
\[o'_A=Q\iota_A,\]
for some $Q$ also to be  found. Next we calculate
\[{\bp}_{,a}=({\by}_{,a})'=4Po_{(A}\iota_B\iota_C\iota_{D)},\]
and continue in the same way to express ${\bq}_{,a},{\br}_{,a}$ and ${\bs}_{,a}$ in the spinor dyad. One more derivative together with (\ref{eq1}) gives a set of identities which fix $P$ and $Q$. 

Defining the metric $\tilde{g}_{ab}$ from $\epsilon_{AB}$ in the standard way as
\[\tilde{g}_{ap}=\tilde{g}_{ABCD.PQRS}=\epsilon_{(A}^{\;\;W} \epsilon_{B}^{\;\;X}\epsilon_{C}^{\;\;Y} \epsilon_{D)}^{\;\;Z}\epsilon_{PW} \epsilon_{QX} \epsilon_{RY}  \epsilon_{SZ} 
\]
we may express all coordinate gradients $X^{\bf{a}}_{\;,a}$ in the spinor dyad and hence calculate all inner products $\tilde{g}^{ab}X^{\bf{a}}_{\;,a}X^{\bf{b}}_{\;,b}$, and then deduce the metric in these coordinates. 
For the covariant metric we find
  \[
(\tilde{g}_{{\bf ab}})=\left(\begin{array}{ccccc}
           \frac{\br^2\bs}{24\bq^5}-\frac{5\br^4}{162\bq^6}-\frac{\bs^2}{72\bq^4} & \frac{\br\bs}{72\bq^4}-\frac{\br^3}{54\bq^5} & \frac{13}{72}\frac{\br^2}{\bq^4}-\frac{\bs}{12\bq^3} & -\frac{\br}{8\bq^3} & \frac{1}{24\bq^2)}\\
          * & \frac{\bs}{24\bq^3}-\frac{\br^2}{18\bq^4} & \frac{\br}{24\bq^3} & -\frac{1}{24\bq^2)} & 0\\
           * & * & -\frac{1}{24\bq^2)} & 0 & 0 \\
           * &* & * &0 & 0\\
           * & *  &*&*&0\\
\end{array}\right),
\]
where asterisked terms are determined by symmetry. It is straightforward to check that this metric is Einstein with $R=-60$ and that it is indeed the metric of \S\ref{sub22} multiplied by the constant factor 1/64. The corresponding 
$\tilde{G}_a^{\;\;bc}$ with indices arranged like that is readily found to be $8G_a^{\;\;bc}$, again with $G_a^{\;\;bc}$ as in \S\ref{subso3}. These changes modify the eigenvalues in (3.16) so that those equations 
become
\be\label{eq22}\tilde{\Delta}F=-\frac{16}{3}F,\;\;\;\tilde{\Box}F=\frac{1}{3}dF,\ee
and our aim is to deduce these equations from the integral expression (\ref{eq23}) for $F$, which can now be written
\[F(X^{\bf{a}})=\oint \tf(x,Z(x,X^{\bf{a}})){\bq}^{1/3}dx.\]

Calculating in abstract indices we obtain
\[\nabla_aF=\oint \left(\frac{\partial \tf}{\partial y}Z_a{\bq}^{1/3}+\frac{1}{3}\tf{\bq}^{-2/3}{\bq}_a\right)dx\]
\[\tilde{\nabla}_a\nabla_bF=\oint\left(\frac{\partial^2 \tf}{\partial y^2}Z_aZ_b+\frac{\partial \tf}{\partial y}\tilde{\nabla}_aZ_b\right){\bq}^{1/3}dx\]
\[+\oint\left(\frac{1}{3}{\bq}^{-2/3}\frac{\partial \tf}{\partial y}(Z_a{\bq}_b+Z_b{\bq}_a)+\frac{1}{3}\tf({\bq}^{-2/3}\tilde{\nabla}_a{\bq}_b-\frac{2}{3}{\bq}^{-5/3}{\bq}_a{\bq}_b)\right)dx.\]
For $\tilde{\Delta} F$ we need to note
\[\tilde{g}^{ab}Z_aZ_b=0=\tilde{\Delta} Z=\tilde{g}^{ab}Z_a{\bq}_b=\tilde{\Delta} {\bq},\;\;\tilde{g}^{ab}{\bq}_a{\bq}_b-24{\bq}^2=0,\]
all of which follow from the expressions for $\tilde{g}^{ab}X^{\bf{a}}_{\;,a}X^{\bf{b}}_{\;,b}$, for then
\[\tilde{\Delta} F=-\frac{2}{9}\oint 24{\bq}^{1/3}\tf dx=-\frac{16}{3}F,\]
as required.

For $\tilde{\Box}F$ we need expressions for $\tilde{G}_a^{\;\;bc}X^{\bf{b}}_{\;b}X^{\bf{c}}_{\;c}$ which can be obtained from the coordinate gradients $X^{\bf{a}}_{\;,a}$ in the spinor dyad. In particular we find
\[\tilde{G}_a^{\;\;bc}Z_bZ_c=0=\tilde{G}_a^{\;\;bc}{\bq}_bZ_c-2{\bq}Z_a=  \tilde{G}_a^{\;\;bc}\tilde{\nabla}_bZ_c+Z_a= \tilde{G}_a^{\;\;bc}\tilde{\nabla}_b{\bq}_c+{\bq}_a=\tilde{G}_a^{\;\;bc}{\bq}_b{\bq}_c+2{\bq\bq}_a.\]
Then
\[\tilde{\Box}F=\int\left(-{\bq}^{1/3}\tf_ydZ+\frac{4}{3}{\bq}^{1/3}{\tf}_ydZ-\frac{1}{3}{\bq}^{-2/3}\tf d{\bq}-\frac{2}{9}{\bq}^{-5/3}\tf(-2{\bq}d{\bq})\right)dx\]
\[=\int\left(\frac{1}{3}{\bq}^{1/3}\tf_ydZ+\frac{1}{9}\tf{\bq}^{-2/3}d{\bq}\right)dx=\frac{1}{3}dF,\]
again as required, establishing the claim.
\koniec
In the Appendix (formula (\ref{big_operator}))
we give explicit expressions for the second order operator
anihilating the functions in the range of (\ref{contour_int}).


The next proposition shows that the two sets of equations
(\ref{moraru_eq}) are not independent:
$\Delta_g F=-\frac{1}{12}F$ is implied by $\Box F=\frac{1}{24}dF$.
We shall establish a slightly more general result applicable to other forms
of integrable $SO(3)$--structures with arbitrary constant and non--zero Ricci scalar
\begin{prop}
\label{theo_tod_2}
Let $(M, g, G)$ be a five--dimensional integrable $SO(3)$ structure such that the Ricci scalar $R$ of
$g$ is constant, and different than zero. Let $F:M\rightarrow \R$ satisfy
\be
G_a^{\;\;bc}\nabla_b\nabla_cF=\lambd\nabla_aF \label{G5}
\ee
for some constant $\kappa$.
Then
\be
\Delta_g F=\mu F \label{G6}, 
\quad\mbox{where}\quad \mu=6\lambd^2+\frac{R}{10}.
\ee
\end{prop} 
\noindent
{\bf Proof.}
Consider the $SO(3)$ structure (\ref{so3structure}), and
trace (\ref{so3structure}) to obtain
\[G_{efa}G^{ef}_{\;\;\;\;b}=\frac{7}{12}g_{ab}\quad
\mbox{and}\quad G_{abc}G^{abc}=\frac{35}{12}.\]
Commute derivatives on $G_{abc}$ to obtain
\be\label{G1}R_{abc}^{\;\;\;\;\;\;(d}G^{ef)c}=0,\ee
where ${R_{abc}}^d$ is the Riemann tensor of $g$.
Define
\[\chi_{abcd}=6G^e_{\;\;ab}G_{cde},\;\;F_{bcad}=\chi_{a[bc]d},\]
and claim
\be\label{G2}\chi_{abcd}=\chi_{(abcd)}+\frac{2}{3}F_{bcad}+\frac{2}{3}F_{bdac},\ee
with
\[\chi_{(abcd)}=6G^e_{\;\;(ab}G_{cd)e}=g_{a(b}g_{cd)}.\]
Expand (\ref{G1}):
\[R_{abc}^{\;\;\;\;\;\;d}G^{efc}+R_{abc}^{\;\;\;\;\;\;e}G^{fdc}+R_{abc}^{\;\;\;\;\;\;f}G^{dec}=0\]
and contract with $G_{efp}$ to deduce
\be\label{G4}R_{abcd}F^{cd}_{\;\;\;\;\;\;pq}=\frac{7}{4}R_{abpq},\ee
after relabeling of indices.

The system of interest is (\ref{G5}, \ref{G6})
Compress notation by writing $F_a=\nabla_aF$ then from (\ref{G5})
\[6\lambd^2F^a=6G^{abc}G_b^{\;\;de}\nabla_c\nabla_dF_e=\chi^{acde}\nabla_c\nabla_dF_e\]
\[=(g^{a(c}g^{de)}+\frac{2}{3}F^{cdae}+\frac{2}{3}F^{cead})\nabla_c\nabla_dF_e\]
Here the first term is
\[\frac{1}{3}(g^{ac}g^{de}+g^{ad}g^{ec}+g^{ae}g^{cd})\nabla_c\nabla_dF_e\]
\[=\frac{1}{3}(\nabla^a\Delta F+2\nabla_c\nabla^aF^c)\]
\[=\frac{1}{3}(\nabla^a\Delta F+2R^{ab}F_b+2\nabla^a\Delta F)\]
\[=\nabla^a(\Delta F+\frac{2}{15}RF).\]
The other two terms become
\[\frac{4}{3}F^{cdae}\nabla_c\nabla_dF_e=-\frac{2}{3}F^{cdae}R_{cdfe}F^f\]
\[=-\frac{2}{3}.\frac{7}{4}R^{ae}_{\;\;\;\;\;fe}F^f=-\frac{7}{6}R^a_{\;\;\;f}F^f=-\frac{7}{30}RF^a.\]
Putting them together
\[6\lambd^2F_a=\nabla_a(\Delta F-\frac{1}{10}RF),\]
whence
\be\label{G7}\mu=6\lambd^2+\frac{R}{10}.\ee
Conversely, a solution $F$ of (\ref{G5}) with some $\lambd$ will necessarily satisfy (\ref{G6}) with the value of $\mu$ given by (\ref{G7}), possibly after adding a constant to $F$.
\koniec
\subsection{Examples}
The general solution of the system (\ref{moraru_eq}) is given by integrating
the cohomology classes of functions on suitable $U\subset\CP^2$
 along conics. An explicit formula
for $F$ seems to be out of reach, but there is a class of solutions of the form
\be
\label{F1}
F=(\gamma_1+\gamma_2p+\gamma_3 r +
( \gamma_4+\gamma_5p)(a+2b))e^{-a-b}+
(\gamma_6+\gamma_7 q+\gamma_8r)e^{-2a},
\ee
where $(\gamma_1, \dots, \gamma_8)$ are arbitrary constants. 
In particular all solutions which are independent of $(q, r)$ are of this form with
$\gamma_3=\gamma_7=\gamma_8=0$. 

 Another class is obtained by looking for $F$ which does not depend on $r$.
Setting $u=\sqrt{p^2e^{2(b-a)}+e^{-2(2a+b)}}$ leads to a general solution of the form
\[
F=\frac{1}{u}e^{-2a}K(u, q)+F_1(a, b, p), \quad \mbox{where} \quad
\frac{\p^2 K}{\p q^2}+\frac{\p^2 K}{\p u^2}=0
\]
and $F_1(a, b, p)$ is of the form (\ref{F1}) with 
$\gamma_3=\gamma_7=\gamma_8=0$.
More examples arise from utilising the $SL(3)$ action on
solutions to (\ref{moraru_eq})
induced by the action (\ref{sl_action}). For example by taking
$
F=re^{-2a}=\frac{A_{23}}{(A_{33})^2}, 
$
and replacing $A$ by $\hat{A}=NAN^T$.

\section{Anti--self--dual conformal structures in dimension four}
\label{sec_asd}
This is the main section of our paper. Given a function $F$ in the range
of the Penrose--Radon transform from \S\ref{sec_radon}, we shall construct
an anti--self--dual metric $\gamma$ on the hypersurface $X$
given by the zero set of $F$ in $M$. The twistor space of $(X, \gamma)$
fibers holomorphically over $\CP^2$, and the twistor curves
project to the four--parameter family of conics in $\CP^2$ such that
the cohomology class corresponding to $F$ vanishes on this family. 
The resulting conformal structure is equivalent to that constructed by
Moraru \cite{moraru},
but our procedure is different and leads to an explicit metric
(which we call barycenter in \S\ref{bary_section}), which admits three
linearly independent solutions to the conformal Killing--Yano
equations. Thus the barycenter metric is conformal to K\"ahler in three different ways. In \S\ref{sec_asd_ex} we shall give several examples of this construction.
We shall prove
\begin{theo1}
\label{theo1_intro}
Let $F:M\rightarrow \R$, and let $X$ be a four--manifold  defined by
\[
X=\{m\in M, F(m)=0\}.
\]
Let $[\gamma]$ be a conformal structure on $X$ such that $V\in \C\otimes TM$ is null
if $\rho(V)=h\otimes\iota$ where $h\in \mbox{Sym}^3(\C^2)$ is such that $<h, \rho(dF)>_3=0$  and 
$\iota\in \C$. Then
\begin{enumerate}
\item The Weyl tensor of $[\gamma]$ is anti-self-dual if $F$ satisfies (\ref{theoremb_eq}).
\item There exists a basis $\{\Omega^{(1)}, \Omega^{(2)}, \Omega^{(3)}\}$ of ${\Lambda^2}_+(X)$
such that $d\Omega^{(i)}=0, i=1, 2, 3$ if and only  if $F$ satisfies (\ref{theoremb_eq}).
\item  The conformal class $[\gamma]$ is tri--K\"ahler:
given any metric $\gamma \in [\gamma]$ there exist three scalar--flat K\"ahler structures
$(\gamma^{(i)}, \Omega^{(i)})$, where $\gamma^{(i)}=|\Omega^{(i)}|_{\gamma}\gamma$.
\end{enumerate}
\end{theo1}
\noindent
This theorem will follow from  Theorems \ref{theo_1_forms} and 
\ref{theo_2_forms}.

Define a four--manifold by
\be
\label{4_mfd}
X=\{m\in M, F(m)=0\},
\ee
where the function $F:M\rightarrow \R$ satisfies (\ref{moraru_eq}).
Let $\rho(dF)\in\mbox{Sym}^4(\C^2)$ be the quartic corresponding
to $dF\in\Lambda^1(M)$ by the $GL(2)$ structure (\ref{para_con_exp}).
Let $H\in\mbox{Sym}^3(\C^2)$ be a cubic and let
$l\in\C^2$ be a linear form.
The conformal structure $[\gamma]$ on $X$ is determined by specifying the null cone to be the set of quartics
\be
\label{null_cone}
{\mathcal{N}}=\{H\otimes l, \quad \mbox{where} \quad <\rho(dF), H>_3=0\}.
\ee
This will be non--degenerate iff the $J$--invariant
of the quartic $\rho(dF)$ given by 
\be
\label{j_invariant}
J=<<\rho(dF), \rho(dF)>_2, \rho(dF)>_4
\ee
does not vanish, as then the space of solutions
to the linear system of equations (\ref{null_cone}) for the four components
of $H$ is two--dimensional.

We shall now rephrase this in the spinor notation (see Appendix A).
Let $F:M\rightarrow \R$ and let $F_{ABCD}$ be defined by
$dF=F_{ABCD}e^{ABCD}$, where $e^{ABCD}$ is defined by (\ref{abcd}).

 A vector $V^{ABCD}$ is 
null on $X$ iff $ V^{ABCD}=h^{(ABC}\iota^{D)}$, where $F_{ABCD}h^{ABC}=0$.
Instead of solving this system for $h^{ABC}$ we use projections
\[
h^{ABC}\rightarrow |F|^2h^{ABC}-2h^{PQR}F_{PQRS}F^{SABC},
\]
of (for example) $o^Ao^Bo^C$ and $\iota^A\iota^B\iota^C$.
To construct the explicit frame for $[\gamma]$ set
\[
H=h_0s^3+3h_1s^2t+3h_2st^2+h_3t^3,
\]
and solve (\ref{null_cone}) for $(h_1, h_2)$ in terms of 
$(h_0, h_3)$. Let $H_0$ be the cubic corresponding to
$(h_0, h_3)=(1, 0)$, and let $H_1$ corresponds to
$(h_0, h_3)=(0, 1)$. Let $V^{00}, V^{10}, V^{01}, V^{11}$
be four vector fields on $M$ corresponding to the quartics
$
H_0s, \, H_0t, \, H_1s, H_1t.
$
By construction,  these vector fields annihilate the one--form $dF$,
and so they span $TX$. To construct the conformal stucture
explicitly, let $V_F\in TM$ be a vector field
such that $V_F\hook S=dF$, where $S$ is the $GL(2)$ structure 
(\ref{para_S}), and let 
$
e^{00}, e^{01},  e^{10}, e^{11}, dV
$
be a basis of $T^*M$ dual to $V^{00}, V^{01}, V^{10}, V^{11}, V_F$.
Then the ASD conformal structure is
\be
\label{The_metric_4d}
[\gamma]=\Omega^2(e^{01}e^{10}-e^{00}e^{11})
\ee
where $\Omega:X\rightarrow \R^+$.
\subsection{Conformal structure from self--dual two--forms}
First recall the `usual' twistor picture \cite{penrose}, but with primed and unprimed indices 
swapped round. Let $(X, \gamma)$ be an oriented Riemannian four--manifold with volume form $\mbox{vol}_X$.  
The integrable twistor distribution is
$L_{A'}=\pi^A\nabla_{A'A}$, and the self--dual (SD) two forms are
$\Sigma^{AB}=1/2 \epsilon_{AB} e^{A'A}\wedge e^{B'B}$. The relation between
$L_{A'}$ and $\Sigma^{AB}$ is
\be
\label{forms_from_lax}
\mbox{vol}_X(L_{A'}, L_{B'}, \cdot, \cdot)=
\epsilon_{A'B'}\pi_{A}\pi_{B}\Sigma^{AB}.
\ee
Conversely, given $\Sigma=\pi_{A}\pi_{B}\Sigma^{AB}$, the twistor distribution
(and so the conformal structure) arises as the kernel of $\Sigma$, and the conformal structure can be 
recovered from the Urbantke formula \cite{urbantke}: 
Let $[\gamma]$ be a conformal structure on a four--manifold $X$ such that
the two--forms $\Sigma^{AB}$ are self--dual. Then $\gamma(V, V)=0$ for
any $\gamma\in [\gamma]$ if and only if
\[
\Sigma_{AB}(V, \cdot)\wedge \Sigma^{BC}(V, \cdot)\wedge {\Sigma_C}^A=0.
\]
In Theorem \ref{theo_1_forms} below we shall reproduce (\ref{forms_from_lax}) for the twistor discribution arising on the space of conics.

We make use of  the isomorphisms
\[
\C\otimes T_mM\cong \sym^4(\C^2), \quad \Lambda^2(\sym^4(\C^2))\cong\sym^6(\C^2)\oplus \sym^2(\C^2)
\]
to introduce a basis of $\Lambda^{2}(M)$
\[
\sigma^{AB}=\frac{1}{48}e^{ACDE}\wedge {e^B}_{CDE}, \quad
\sigma^{ABCDEF}=\frac{1}{8}e^{G(ABC}\wedge {e^{DEF)}}_{G},
\]
where $\sigma^{AB}$ and $\sigma^{ABCDEF}$ are defined in the Appendix.

Let us define three two--forms $\Sigma^{AB}=\Sigma^{(AB)}$ on $M$ by
\be
\label{spinor_sigma}
\Sigma^{AB}={F^{AB}}_{CD}\sigma^{CD}+\frac{15}{2}F_{CDEF}\sigma^{ABCDEF}.
\ee
These will pull back to two--forms on 
the four--manifold $X$ defined by (\ref{4_mfd}). 
\begin{lemma}
Let ${\Gamma^{A}}_{B}$ be the spin--connection of $(M, g)$. The following identity
\be
\label{identity}
d\Sigma^{AB}+\frac{1}{2}{\Gamma^{(A}}_{C}\wedge\Sigma^{B)C}=
-\sqrt{6}{\Box^{AB}}_{CD}F\star\sigma^{CD}+5\sqrt{6}\Box_{CDEF}F\star
\sigma^{ABCDEF}-\frac{5\sqrt{6}}{6}\Delta_g F\star \sigma^{AB}
\ee
holds for any $F:M\rightarrow\R$.
\end{lemma}
\noindent
{\bf Proof.} 
We find that
\be
\label{connection_j}
d\sigma^{AB}+\frac{1}{2}{\Gamma^{(A}}_{C}\wedge \sigma^{B)C}=0, \quad
d\sigma^{ABCDEF}+\frac{3}{2}{\Gamma^{(A}}_{G}\wedge \sigma^{BCDEF)G}=0,
\ee
where
\[
{\Gamma^1}_1=-{\Gamma^0}_0=\frac{1}{4}(e^{0000}-e^{1111}), \quad{\Gamma^1}_0=-e^{0001}, \quad {\Gamma^0}_1=-e^{0111}.
\]
In five dimensions $\Lambda^2(M)\cong\Lambda^{3}(M)$ by Hodge isomorphism, and we  take
$\star\sigma^{AB},\star \sigma^{ABCDEF}$ to be the basis of $\Lambda^3(M)$.
The following identities can be established by an explicit computation
\begin{eqnarray}
\label{id_j}
e_{ABCD}\wedge \sigma^{EF}&=&-\frac{\sqrt{6}}{2}\star {\sigma_{ABCD}}^{EF}
-\frac{2\sqrt{6}}{5}{\epsilon^{(E}}_{(A}{\epsilon^{F)}}_{B}\star\sigma_{CD)},\\
e_{ABCD}\wedge \sigma^{EFGHIJ}&=&
-\frac{\sqrt{6}}{5}
{\epsilon^{(E}}_{(A}{\epsilon^F}_{B} {\epsilon^G}_{C}{\epsilon^H}_{D)} 
\star\sigma^{IJ)}+\sqrt{6} 
{\epsilon^{(E}}_{(A}{\epsilon^F}_{B}\star{\sigma^{GHIJ)}}_{CD)}.\nonumber
\end{eqnarray}
The identity (\ref{identity}) can now be verified directly.
\koniec
\subsection{The double fibration picture}
Consider a double fibration
\[
M_\C\longleftarrow {\mathcal F} \stackrel{\rho}\longrightarrow \CP^2,
\]
where ${\mathcal F}\subset M_\C \times \CP^2$ is the six--dimensional
manifold of incident pairs $(m, \xi)$ such that
$\xi\in C_m$, where $\xi\in \CP^2$, and the conic $C_m\subset \CP^2$ corresponds to a point $m\in M_\C$. A point $\xi\in \CP^2$ corresponds to a hyper--surface
$\hat{\xi}\subset TM_\C$ which is totally null in the sense of \cite{DT}:
the normal vector to $\hat{\xi}$ is a polynomial with a quadruple root.
The  map $\rho$ is the quotient by the rank-four distribution 
${\mathcal O}_{\CP^1}(-1)\otimes\C^4$ given by
${\mathcal D}_{ABC}=\pi^{D}\nabla_{ABCD}$. Consider the projection of this distribution to a 
rank-two distribution on ${\mathcal F}$ 
\be
\label{LABC}
{\mathcal D}_{ABC}\rightarrow L_{ABC}=|F|^2{\mathcal D}_{ABC}
-2{F_{ABCD}}F^{DPQR}{\mathcal D}_{PQR}.
\ee
In the next Theorem we shall establish the  integrability of the distribution
${\mathcal D}_{ABC}$ assuming that the system (\ref{moraru_eq}) holds for $F$.
\begin{theo} Let $[\gamma]$ be a conformal structure on the four--manifold $X$
(\ref{4_mfd}) such that $V\in TX$ is null iff $V^{ABCD}=h^{(ABC}\iota^{D)}$, where $F_{ABCD}h^{ABC}=0$. Then $[\gamma]$ is ASD if equations (\ref{moraru_eq}) hold.
\end{theo}
\noindent
{\bf Proof.}
In the proof we shall use the constant rescaling  of the metric by $1/64$ which leads to equations 
(\ref{eq22}).

Given $M$ and the zero-locus $X$ of $F$, consider the two-dimensional distribution spanned by $h^{(ABC}\iota^{D)}$ on $X$ where $\iota^A$ is the spinor field on $M$ introduced in (\ref{special_spinor}) for one fixed choice of $x$ 
and $h^{ABC}$ is any solution of
\[h^{ABC}F_{ABCD}=0,\]
then this distribution is integrable provided $F$ satisfies the system
\be\label{sys1}\Box_{ABCD}F=\lambd F_{ABCD},\;\;\;\Delta F=\mu F,\ee
with $\lambd=1/3$ and $\mu=-16/3$.

 Evidently vectors of the form $h^{(ABC}\iota^{D)}$ are tangent to $X$.
 Choose two such vectors
\[X^a=h^{(ABC}\iota^{D)},\;\;\;Y^a=\tilde{h}^{(ABC}\iota^{D)},\]
and consider their commutator $Z=[X, Y]$, or
\[Z^{ABCD}:=X^{PQRS}\nabla_{PQRS}Y^{ABCD}-Y^{PQRS}\nabla_{PQRS}X^{ABCD},\]
where $\nabla_{ABCD}$ is the metric covariant derivative on $M$. We wish to show that 
$Z^{ABCD}=\hat{h}^{(ABC}\iota^{D)}$ for some $\hat{h}^{ABC}$ annihilating $F_{ABCD}$. From \cite{MC} we quote the formula 
\be\label{2}\iota^S\nabla_{PQRS}\iota_D=\frac12\iota_P\iota_Q\iota_Ro_D +\gamma\iota_P\iota_Q\iota_R\iota_D,\ee
for some $\gamma$ which is known but will turn out to be irrelevant, and $(o^A,\iota^A)$ is a normalised spinor dyad. We calculate
\[Z^{ABCD}=X^{PQRS}\nabla_{PQRS}Y^{ABCD}-Y^{PQRS}\nabla_{PQRS}X^{ABCD}\]
\[=\iota^S(h^{PQR}\nabla_{PQRS}\tilde{h}^{(ABC}-\tilde{h}^{PQR}\nabla_{PQRS}h^{(ABC})\iota^{D)}\]
\[+\iota^S(h^{PQR}\tilde{h}^{(ABC}-\tilde{h}^{PQR}h^{(ABC})\nabla_{PQRS}\iota^{D)}.\]
It will be convenient to introduce a spinor field $\chi^{QR}_{\;\;\;\;\;\;BC}$, symmetric on each pair of spinor indices and with the interchange symmetry, by
\[h^{PQR}\tilde{h}_{ABC}-\tilde{h}^{PQR}h_{ABC}:=\delta^{(P}_{(A}\chi^{QR)}_{\;\;\;\;\;\;BC)}.\]
From the definition of $\chi^{QR}_{\;\;\;\;\;\;BC}$ it follows that 
\[0=F_{PQRS}(h^{PQR}\tilde{h}_{ABC}-\tilde{h}^{PQR}h_{ABC})=F_{PQRS}\delta^{(P}_{(A}\chi^{QR)}_{\;\;\;\;\;\;BC)}=F_{QRS(A}\chi^{QR}_{\;\;\;\;\;\;BC)},\]
which also implies that
\[F_{QRAB}\chi^{QR}_{\;\;\;\;\;\;CD}=\frac12\phi(\epsilon_{AC}\epsilon_{BD}+\epsilon_{BC}\epsilon_{AD})\]
for some $\phi$, when by tracing on $AC$
\[F_{QRSB}\chi^{QRSD}=\frac{3}{2}\phi\delta_B^D,\]
and, by tracing again,
\[\phi=\frac{1}{3}F_{PQRS}\chi^{PQRS}.\]
Thus
\be\label{0}F_{QRAB}\chi^{QR}_{\;\;\;\;\;\;CD}=\frac{1}{6}(F_{PQRS}\chi^{PQRS})(\epsilon_{AC}\epsilon_{BD}+\epsilon_{BC}\epsilon_{AD}).\ee
With these we deduce that 
\[Z^{ABCD}=\phi^{(ABC}\iota^{D)} \]
with
\be\label{9}
\phi^{ABC}=\iota^S(h^{PQR}\nabla_{PQRS}\tilde{h}^{ABC}-\tilde{h}^{PQR}\nabla_{PQRS}h^{ABC})+\iota_Q\iota_R\chi^{QR(AB}(\frac12o^{C)}+\gamma\iota^{C)}),
\ee
where the last term uses (\ref{2}). It remains to check that $\phi^{ABC}F_{ABCD}$ vanishes. From (\ref{9}) we calculate
\begin{eqnarray*}
&&F_{ABCD}\phi^{ABC}=\\
&&F_{ABCD}\left(\iota^S(h^{PQR}\nabla_{PQRS}\tilde{h}^{ABC}-\tilde{h}^{PQR}\nabla_{PQRS}h^{ABC})+\iota_Q\iota_R\chi^{QR(AB}(\frac12o^{C)}+\gamma\iota^{C)}\right)
\end{eqnarray*}
and then, using $h^{ABC}F_{ABCD}=0=\tilde{h}^{ABC}F_{ABCD}$, 
\begin{eqnarray}
\label{sys2}
&&F_{ABCD}\phi^{ABC}=\\
&&-\iota^S(h^{PQR}\tilde{h}^{ABC}-\tilde{h}^{PQR}h^{ABC})\nabla_{PQRS}F_{ABCD}+F_{ABCD}\iota_Q\iota_R\chi^{QR(AB}(\frac12o^{C)}+\gamma\iota^{C)})\nonumber.
\end{eqnarray}
Now we need to notice, again by symmetry, that there is a spinor field $Q^{RS}_{\;\;\;\;\;\;CD}$ symmetric in both pairs of spinor indices and with the interchange symmetry satisfying
\[\nabla^{AQRS}F_{ABCD}=\nabla^{AQRS}\nabla_{ABCD}F=\delta^{(Q}_{(B}Q^{RS)}_{\;\;\;\;\;\;CD)}.\]
Contractions then show that
\[Q^{RS}_{\;\;\;\;\;\;CD}=\frac{3}{2}\Box^{RS}_{\;\;\;\;\;\;CD}F+\frac{1}{6}\Delta F(\delta^R_C\delta^S_D+\delta^R_D\delta^S_C).\]
when with the aid of the system (\ref{sys1}),
\[Q^{RS}_{\;\;\;\;\;\;CD}=\frac{3\lambd}{2}F^{RS}_{\;\;\;\;\;\;CD}+\frac{\mu}{3}\delta^{(R}_{(C}\delta^{S)}_{D)}F.\]
With these in hand, (\ref{sys2}) reduces to
\[F_{ABCD}\phi^{ABC}=\frac{1}{2}\lambd(\chi_{ABCD}F^{ABCD})\iota_D-\frac{1}{6}(\chi_{ABCD}F^{ABCD})\iota_D,\]
where we have dropped the term proportional to $\mu$ as it vanishes on $X$. Evidently this vanishes if $\lambd=1/3$ (and we have seen elsewhere that with this value of $\lambd$ the other part of the system 
(\ref{sys1}) automatically holds with $\mu=-16/3$).
\koniec

The leaves of the distribution (\ref{LABC}) are $\alpha$--surfaces
of $(X, [\gamma])$.
The set of conics through a given point in
$\CP^2$ in a  given direction is  a three--dimensional surface in $M$. The intersection of this
surface with the hyper--surface $F=0$ is two--dimensional. This is an $\alpha$--surface
in $X$ or equivalently a point in the twistor space 
${\mathcal Z}$ of  $(X, [\gamma])$.

In the next Theorem we shall  give a direct way to construct the ASD conformal structure on $X$ in terms of a preferred basis of self--dual two--forms on $X$.
Let $\pi:X\rightarrow M$ be the map given by (\ref{4_mfd}).
\begin{theo}
\label{theo_1_forms}
Let $[\gamma]$ be the conformal structure on the four--manifold $X$
(\ref{4_mfd}) such that
$V\in \C\otimes TX$ is null iff $V^{ABCD}=h^{(ABC}\iota^{D)}$, where $F_{ABCD}h^{ABC}=0$. Then the two--forms
$\Sigma^{AB}$ given by (\ref{spinor_sigma}) pull back to  two--forms which are self--dual w.r.t the orientation given by
the pull back of $\star_5 dF$ from $M$ to $X$.

Conversely, let $\Sigma^{AB}$ be given by (\ref{spinor_sigma}). Then
\[
\pi^*(\Sigma^{(AB}\wedge\Sigma^{CD)})=0,
\] 
and so 
$\Sigma^{AB}=(1/2)\epsilon_{A'B'}e^{AA'}\wedge e^{BB'}$ for some tetrad
$e^{AA'}$ on $X$ which is non--degenrate iff ${J}\neq 0$. The corresponding conformal structure is given by $[\gamma]$.
\end{theo}
\noindent
{\bf Proof.}
Consider  the rank 2 distribution  $L_{ABC}$ given by (\ref{LABC}).
Therefore for any $v^{ABC}$ the vector field $V=v^{ABC}L_{ABC}$ is a null vector
field on $(X, [\gamma])$ in the sense of Moraru:
\[
V^{ABCD}=\pi^{(A}h^{BCD)}, \quad\mbox{where}\quad 
h^{BCD}= |F|^2v^{BCD}-2F^{BCDE}{F}_{EPQR}\;v^{PQR}
\]
so that $h^{BCD}F_{ABCD}=0$, where we have used the identity
$2F_{ABCD}F^{ABCE}= {\delta_D}^E|F|^2$.

Let $\mbox{vol}_X=\pi^*(\star_5 dF)$ be a volume form on $X$. Consider
the ${\mathcal O}_{\CP^1}(2)$--valued two form on $\CP^1\times X$ given by
\be
\label{s_form}
{\mathcal S}=(\star_5dF)(L_{000}, L_{111}, \cdot, \cdot).
\ee
We have verified using MAPLE 
that
${\mathcal S}=s \pi^*(\pi_A\pi_B\Sigma^{AB})$, where the scalar multiple $s$ is, up to a constant numerical factor,  given by
\[
s=
\pi^*(|F|^2 F_{PQRS}{F^{PQ}}_{MN} o^R o^M \iota^S\iota^N),
\]
and the two--forms $\Sigma^{AB}$ are given by (\ref{spinor_sigma}).
To establish the second part compute
\[
\Sigma^{00}\wedge \Sigma^{00}\wedge  dF=\Sigma^{00}\wedge \Sigma^{01}\wedge dF=
\Sigma^{11}\wedge \Sigma^{11}\wedge dF=\Sigma^{11}\wedge \Sigma^{01}\wedge dF=0,
\]
\[
\star(\Sigma^{00}\wedge \Sigma^{11}\wedge dF)=
-2\star(\Sigma^{01}\wedge \Sigma^{01}\wedge dF)=\frac{25\sqrt{6}}{2592}
{J}.
\]
\koniec
\subsection{The barycenter metric, and the tri--K\"ahler structure}
\label{bary_section}
According to Moraru \cite{moraru} there is a sphere of  scalar--flat--K\"ahler 
metrics in every conformal class arising from conics. 
In Theorem \ref{theo_2_forms} below we shall construct these K\"ahler forms explicitly, but
let us first explain why one should expect them to exist from the twistor perspective.

Pick a section of 
$\OO_{\CP^2}(1)$ which is of the form $\omega=W_i Z^i$ for some 
$[W]\in {\CP^2}^*$, and restrict it to a conic $C_m$ given by 
(\ref{parametrisation}) 
and such that $F(m)=0$. This gives a quadratic polynomial - a section of
$\OO(2)$ restricted to a twistor line - given by $\omega_{AB}\pi^A\pi^B$,
where $\omega_{AB}$ depends on $m\in M$, as well as $W$. There is a two--parameter 
family of $W$s, as linear functions on $\CP^2$ are defined up to scale
of $\C^3$, so there is a (at least) three--dimensional space of solutions
to the twistor equation $\nabla_{A'(A}\omega_{BC)}=0$ on the 
four--manifold $X$ given by (\ref{4_mfd}).

Let $Z=[Z^1, Z^2, Z^3]$ be homogeneous coordinates on $\CP^2$, and let
a rational parametrisation of the conic $ZAZ^T=0$ be 
${Z}^{i}=Z^{i}_{AB}\pi^A\pi^B$, where $i=1, 2, 3$, and $Z^{i}_{AB}$
are functions of the components of the symmetric determinant--one matrix ${A}$
which defines the conic. Set
\[
\theta^{(i)}=(Z^{(i)}_{AB}Z^{(i)AB})^{-3/2}\quad\mbox{(no summation over $i$}).
\]
\begin{theo}
\label{theo_2_forms}
The two--forms
\be
\label{three_2_forms}
\Omega^{(i)}=\theta^{(i)}Z^{(i)}_{AB}\Sigma^{AB}\quad\mbox{(no summation over $i$})
\ee
pull--back to self--dual, closed two--forms on $X$ if and only if 
$F$ satisfies (\ref{moraru_eq}).
\end{theo}
\noindent
{\bf Proof.}
Let $\pi:X\rightarrow M$ be given by
$F=0$, where $F:M\rightarrow \R$. Then, for any differential form $\Omega$, 
$\pi^*\Omega=0$ iff $(\Omega\wedge dF)|_{F=0}=0$.
Consider the parametrisation  (with $\pi_{A}=[s, t]$)
\[
Z=
[e^{a+b}(s^2-t^2), 2e^{-b}st-e^{a+b}p(s^2-t^2),  e^{a+b}(pr-q)(s^2-t^2)-2re^{-b}st
+ie^{-a}(s^2+t^2)],
\]
and compute 
\[
\theta^{(1)}=e^{-3a-3b}, \quad \theta^{(2)}=(4p^2e^{2a+2b}+4e^{-2b})^{-3/2},\quad
\theta^{(3)}=(2e^{2a+2b}(pr-q)^2+2e^{-2a}+2r^2e^{-2b})^{-3/2}.
\]
Now use the identity (\ref{identity}) to compute 
\[
(d\Omega^{(i)}\wedge dF)|_{F=0},
\]
and find, using MAPLE, that it vanishes for $i=1, 2, 3$ if $\Box F=(1/24)dF$ 
(this is the second set of equations (\ref{moraru_eq}))
and $\Delta_g F=c F$ for any $c$. Thus the constant $c$ (which should be $-1/12$ to agree with (\ref{moraru_eq})) has not been determined.
We have however shown in Proposition \ref{theo_tod_2} that this constant is determined by the second set of equations together with the scalar 
curvature of $g$.
\koniec

Let $\gamma\subset [\gamma]$
be any metric in the conformal class. Then the conformal class $[\gamma]$ contains
three scalar--flat K\"ahler structures $(\gamma^{(i)}, \Omega^{(i)})$, where
\be
\label{trikah}
\gamma^{(i)}=|\Omega^{(i)}|_\gamma\; \gamma, \quad\mbox{where} \quad|\Omega|_\gamma=
\sqrt{\Omega_{\alpha\beta}\Omega_{\gamma\delta}\gamma^{\alpha\gamma}\gamma^{\beta\delta}}.
\ee
For any metric $\gamma$ in the ASD conformal class $[\gamma]$ 
construct the {\em barycenter} metric
\be
\label{bary_center}
{\gamma}_B=
(|\Omega^{(1)}|_\gamma \cdot|\Omega^{(2)}|_\gamma \cdot
|\Omega^{(3)}|_\gamma )^{1/3}\gamma.
\ee
Then
\begin{eqnarray*}
{\Omega_B}^{(3)}&=&
\frac{\Big(|\Omega^{(1)}|_\gamma \cdot|\Omega^{(2)}|_\gamma\Big)^{\frac{1}{2}}}{|\Omega^{(3)}|_\gamma}
\Omega^{(3)},\\
{\Omega_B}^{(2)}&=&
\frac{\Big(|\Omega^{(1)}|_\gamma\cdot |\Omega^{(3)}|_\gamma\Big)^{\frac{1}{2}}}{|\Omega^{(2)}|_\gamma}\Omega^{(2)},\\
{\Omega_B}^{(1)}&=&
\frac{\Big(|\Omega^{(3)}|_\gamma \cdot|\Omega^{(2)}|_\gamma\Big)^{\frac{1}{2}}}{|\Omega^{(1)}|_\gamma}\Omega^{(1)}
\end{eqnarray*}
are self-dual conformal Killing--Yano tensors for ${\gamma_B}$.
Each such form gives rise \cite{DTkahler} to a K\"ahler class, and the resulting tri--K\"ahler structure is (\ref{trikah})
\subsection{Examples of tri--K\"ahler metrics}
\label{sec_asd_ex}
In this Section we shall construct some  examples of 
anti--self--dual conformal structures containing (at least) three K\"ahler
metrics, and
corresponding to particular subcases of the solution (\ref{F1}).
\subsubsection{An ASD Einstein example}
Consider
\[
F=e^{-2a}+pe^{-a-b},
\]
which is (\ref{F1}) with $\gamma_2=\gamma_6=1$, and the remaining constants equal to zero.
The quartic corresponding to $dF$ via the isomorphism
$\C\otimes T^*M=\mbox{Sym}^4({\spp})$ is
\[
(pe^{-a-b}+2e^{-2a})s^4+4e^{-2a-3b}s^3t+(4e^{-2a}
-2pe^{-a-b})s^2t^2-4e^{-2a-3b}st^3+
(pe^{-a-b}+2e^{-2a})t^4.
\]
Computing the $J$--invariant of this quartic, and
restricting it to the surface $F=0$ yields
$
e^{-6(a+b)},
$
which is nowhere zero. Therefore the  
condition (\ref{null_cone}) has two linearly independent solutions $H_0$ and $H_1$ which lead to
a non--degenerate frame $e^{00},e^{01},e^{10},e^{11}$ and the conformal
structure (\ref{The_metric_4d}) on the four--manifold $X$. It is given by 
\begin{eqnarray*}
e^{00}&=&(e^{3b}-e^{-3b})da-(e^{3b}+2e^{-3b})db
-i(e^{2a+4b}dq+(1+2e^{6b})e^{a-b}dr),\quad
e^{11}=-\overline{e^{00}}\\
e^{10}&=&2da+db+i(e^{2a+b}dq+e^{a+2b}dr),\quad
e^{01}=\overline{e^{10}},
\end{eqnarray*}
and it
gives rise to a metric (\ref{The_metric_4d}).
It is possible to chose the conformal factor 
$\Omega$ such that
the resulting metric is ASD and Einstein, with scalar curvature equal to $-24$
\be
\label{einstein}
\gamma=(da+2db)^2+e^{6b}(da-db)^2+e^{4a+8b}dq^2+4e^{3a+9b}dqdr
+(4e^{6b}+1)e^{2a+4b}dr^2.
\ee
This metric is K\"ahler, but with opposite orientation:
the ASD K\"ahler two form is
\[
\Sigma=e^{2a+4b}(da+2db)\wedge dq+e^{5b+a}(da+5db)\wedge dr.
\]
All anti--self--dual Einstein manifolds which are K\"ahler with opposite 
orientation have constant holomorphic sectional curvature, i.e. 
they are diffeomorphic to $\CP^2$ with its Fubini--Study metric,
its non--compact form $\widetilde{\CP}^2=
SU(2, 1)/U(2)$ with the Bergman metric, or flat space $\C^2$.
Our metric has negative scalar curvature so it is the Bergman space  
$\widetilde{\CP}^2$.  Setting $y=e^{b-a}, z=e^{-a-2b}$ puts it in the form 
\be
\label{hyperb}
\gamma=\frac{1}{z^2}\Big({z^2}h_{\HH^3}+\frac{1}{z^{2}}(dq+2ydr)^2\Big), 
\quad\mbox{where}\quad h_{\HH^3}=\frac{dr^2+dy^2+dz^2}{z^2},
\ee
and the 8--dimensional group of isometries can be constructed explicitly. 
Any ASD Einstein metric with symmetry is conformal to a K\"ahler metric
\cite{DTkahler}, so (\ref{einstein}) contains eight K\"ahler metrics with 
SD K\"ahler forms
in its 
conformal class. Only three of these correspond to the
tri--K\"ahler structure arising from conics. Below we shall examine one of these three, and put it in the canonical $SU(\infty)$--Toda form.
The conformaly rescaled metric 
\[
\tilde{\gamma}=e^{-4a-8b}\gamma
\]
has vanishing scalar curvature, and is K\"ahler with
the `correct' orientation. 
The SD K\"ahler form is
\[
\Omega=e^{-2a-4b}(da+2db)\wedge dq+3e^{-3a-3b}(da+db)\wedge dr.
\]
To recognise this metric we shall put it into 
a general framework of \cite{L91}. 
Any scalar--flat K\"ahler metric with symmetry can be locally put in the form
\be
\label{toda_m}
\tilde{\gamma}={P}(e^u(dr^2+dy^2)+d\zeta^2)+\frac{1}{P}(d q +\alpha)^2
\ee
where $u=u(r, y, \zeta)$  is a solution of the $SU(\infty)$ Toda equation
\[
u_{rr}+u_{yy}+({e^u})_{\zeta\zeta}=0,
\]
the function $P$ is a solution to the linearised $SU(\infty)$ Toda, and
$\alpha$ is a one--form which satisfies the generalised monopole equation.
Set $z^2=2\zeta$, and consider $\tilde{\gamma}=4\zeta^{2}\gamma$, where $\gamma$ is given by
(\ref{hyperb}). Then $\tilde{\gamma}$ is of the form (\ref{toda_m})
where $e^u=2\zeta$, and $P=1$.
\subsubsection{Towards the flat model.} Consider the one--parameter family of solutions
\[
F=e^{-2a}+(p+\kappa r)e^{-a-b}
\]
which is a special case of (\ref{F1}), and 
which contains the conformal class (\ref{einstein})
as the particular case $\kappa=0$. Using the conformal factor
$
\Omega^2=
6^{12}({e^{6b}(\kappa^2+1)+1})e^{-8a-6b}
$
gives an ASD Einstein metric with scalar curvature $-24(\kappa^2+1)$.
It is K\"ahler with the opposite orientation, and the ASD K\"ahler form
\[
e^{2a+4b}(1+\kappa r)(da+2db)\wedge dq+e^{5b+a}(da+5db)\wedge dr
-3\kappa e^{6b} da\wedge db.
\]
The metric $\hat{\gamma}=e^{-4a-8b}\gamma$ is scalar--flat and K\"ahler. The analytic continuation 
$(q, r)\rightarrow(iq, ir)$
of this example to imaginary $\gamma$ gives a one parameter family of $SL(3)$--invariant
Einstein metrics in neutral signature on $SL(3, \R)/GL(2, \R)$. 
The special case $\kappa=i$ gives a flat metric in neutral signature.
\subsubsection{An ASD Ricci--flat example.} The most general element of
$\mbox{Ker}(\hat{\Box})\cap\mbox{Ker}(\triangle+(1/12)Id)$ which depends only on $(a, b)$ is given 
by a special case of (\ref{F1})
\[
F=(\gamma_1+\gamma_2(a+2b))e^{-a-b}+\gamma_3 e^{-2a},
\]
where $\gamma_1, \gamma_2, \gamma_3$ are 
constants,
and $\gamma_2\neq 0$ for $J(dF)\neq 0$. 
This class is characterised by invariance under the three--dimensional group
generated by the Heisenberg algebra $(X_1, X_2, X_3)$, in (\ref{5d_killing}).

Translating 
$(a, b)$ using the $2$--parameter group generated by
$X_4, X_5$ in (\ref{5d_killing}) can be used  
to set $\gamma_1=0$ , and an overall rescaling of $F$ 
allows setting $\gamma_2=1$, which leaves a one--parameter family of solutions depending on one  constant  $\kappa\equiv\gamma_3$ which must be non--zero for $J\neq 0$. The four--manifold (\ref{4_mfd})
is a hypersurface in $M$ parametrised by
\[
a=\frac{1}{3}(s-2\ln{s}+2\ln{\kappa}), \quad
b=\frac{1}{3}(s+\ln{s}-\ln{\kappa}).
\]
The ASD conformal structure (\ref{The_metric_4d}) contains the Ricci--flat metric 
\[
\gamma=(s-1)(e^{-2s}ds^2+\kappa^2e^{-2s}dr^2+ dp^2)+(s-1)^{-1}\kappa^{2}(dq-pdr)^2.
\]
The constant $\kappa$ can be set to $1$ by rescaling $(q, r)$. Setting $(X+iY=e^{-s+ir}, Z=p)$
we recognise this as a Gibbons--Hawking metric
\[
\gamma=V(dX^2+dY^2+dZ^2)+V^{-1}(dq+\Phi)^2,
\]
where $dV=*_3 d\Phi$,  and the harmonic function $V$ on $\R^3$ given by $V=-\ln{(X^2+Y^2)}-1$. Thus we have arrived at the 
Ooguri--Vafa metric \cite{OV} which arises by putting $N$ centres in the 
ALF Gibbons--Hawking gravitational instanton on a line, and taking a limit when $N\rightarrow\infty$.
\subsubsection{A Scalar--flat K\"ahler example}
We shall give one more example, where the barycentre metric 
(\ref{bary_center}) does not contain either Einstein or Ricci-flat metrics
in its conformal class.

Consider 
\[
F=pe^{-a-b}+re^{-2a}.
\]
which is (\ref{F1}) with $\gamma_2=\gamma_8=1$, and remaining constants set to zero. The $J$--invariant of the quartic
$Q(dF)$ is a constant multiple of $
e^{-6a-6b}r.
$
Thus, using $(a, b, r, q)$ as local coordinates on $X$ we expect
the conformal class to degenerate on the hyper--surface $r=0$.

The resulting ASD conformal class (\ref{The_metric_4d}) admits
three scalar--flat K\"ahler metrics with K\"ahler forms
given by (\ref{three_2_forms}). One of these is 
\be
\label{sf_example}
\gamma=\frac{r}{r^2+y^2}(2z(r^2+y^2)(dr^2+dy^2)+dz^2)+
\frac{r^2+y^2}{r}\Big(dq-y^2dy+2rydr+\frac{y}{r^2+y^2}dz\Big)^2,
\ee
where 
\[
y=\frac{b}{a}, \quad z=\frac{1}{2a^2b^4}.
\]
The K\"aher form is given by
\[
\Omega=2d(ryz)\wedge dr+dz\wedge d\Big(q-\frac{y^3}{3}\Big).
\]
The singularity on the $r=0$ surface is a fold in the sense of Hitchin \cite{H2}. The points on this surface correspond to twistor
curves with normal bundle $\OO\oplus\OO(2)$. The 
metric blows up, but
the K\"ahler form remains regular and drops its rank.

Comparing this expression for $\gamma$ with the general form
of the scalar-flat K\"ahler metric with symmetry (\ref{toda_m}) allows us to read--off
the corresponding solutions to the $SU(\infty)$ Toda equation, and
its linearisation:
\[
e^u=2z(r^2+y^2), \quad P=\frac{r}{r^2+y^2}.
\]
\subsubsection{ASD Einstein cohomogeneity--one metrics}
Moraru \cite{moraru} claims (without giving a proof) 
 that all ASD conformal classes which contain
three different K\"ahler (but not hyper--K\"ahler) metrics arises from 
some solution to (\ref{moraru_eq}).

It is known \cite{tod_p6} that all  ASD Einstein cohomogeneity--one metrics arise from
the Painlev\'e VI equation with parameters $(1/8, -1/8, 1/8, 3/8)$.
This case is actually not transcendental, and all solutions can be expressed in terms of the  Weierstrass elliptic function. This is because this particular
PVI is related by a B\"acklund transformation \cite{okamoto}
 to PVI with parameters
$(0, 0, 0, 1/2)$ and that case corresponds to a projectively--flat projective 
structure and has been solved by Picard.

With any symmetry generator we can associate a K\"ahler scale
(by taking the SD derivative), so (if the claim of \cite{moraru} is right) it should arise
from the linear system (\ref{moraru_eq}), and a contour integral formula.
If the Painlev\'e
solution was transcendental it would be a
contradiction, but as it is not, it only shows that
there is some integral formula for the elliptic functions. 
\section{Conics in $\CP^2$ from lines in $\CP^3$}
\label{quad_map_s}

The anti-self-dual tri--K\"ahler metrics arising from Theorem \ref{theo1_intro}
admit twistor spaces  which
holomorphically fiber over $\CP^2$ (this is Moraru's starting point, 
\cite{moraru}), and the twistor curves
with the normal bundle $\OO(1)\oplus\OO(1)$
project to the four--parameter family  conics in $\CP^2$ such that
the cohomology class corresponding to the function 
(\ref{contour_int}) vanishes on this family.
Some twistor spaces with this property appeared in \cite{Cref1}.
In this Section we shall consider a simple case where
the twistor space
$\CP^3\setminus{\mathcal C}$, where ${\mathcal C}\cong\CP^1$ is the rational normal curve.
The holomorphic projection to $\CP^2$ can in this case be described  by classical
projective geometry of Poncelet pairs.

In this section we shall use the notation introduced in \S\ref{sub_transv}
and regard $V_l=\mbox{Sym}^l(\C^2)$ as $(k+1)$--dimensional representation
space of $SL(2, \C)$.

\subsection{A quadratic map}
Consider $\CP^3$ as the projectivisation $\PPP(\mbox{Sym}^3(\C^2))$ of the space $V_3$ of homogeneous cubic polynomials.  We will define a holomorphic projection to $\CP^2=\PPP(\mbox{Sym}^2(\C^2))$ by picking three quadrics in $\CP^3$
and declaring them to be the homogeneous coordinates on $\CP^2$. This triple of quadrics gives a point in $\CP^2$ which (perhaps rather confusingly) is identified with a quadric as $\CP^2=\PPP(V_2)$.

Let ${\mathcal C}$ be the rational normal curve in $\CP^3$ consisting of all
cubics with a triple root $(ax-by)^3\in V_3$. Let the quadratic map $Q$
be defined on the complement of the rational normal curve by
(in this section we shall use the notation from \S\ref{sub_transv},
but applied to homogeneous polynomials in two variables $[x, y]$)
\be
\label{map}
Q:\CP^3\setminus{\mathcal C}\rightarrow \CP^2, \quad
Q(p)=<p, p>_2.
\ee
 Thus in spinor notation (see Appendix), the quadric (a point
in $\CP^2$) corresponding to a cubic $p_{ABC}$ (a point in $\CP^3/{\mathcal C}$)
is ${p_{A}}^{DE}p_{BDE}$. The map (\ref{map}) is given by a choice of
three quadrics in $\CP^3$. The choice is not generic, as the zero 
locus of a generic triple of quadrics in four variables is eight points in $\CP^3$, and our triple vanishes on a curve. 

In Lemmas \ref{lemma1} and \ref{lemma2}
we shall give two more equivalent characterisations of the map (\ref{map}).
In the Lemma \ref{lemma1} below we shall
make use of the map $\CP^1\times \CP^1\rightarrow\CP^2$ which assigns
a quadratic polynomial (up to an overall scale) to a pair of roots, i. e.
$(\alpha_A, \beta_B)\rightarrow \alpha_{(A}\beta_{B)}$. This is a double covering branched over a conic ${\mathcal B}\subset \CP^2$ which is the locus of all points
corresponding to quadratics with a repeated 
root\footnote{Using affine coordinates 
\[
(u, v)\in\CP^1\times\CP^1\rightarrow (z-u)(z-v)\rightarrow
[1, -u-v, uv].
\]
The branch conic $\mathcal B$ corresponding to the diagonal $u=v$ is
${(Z^2)}^2-4Z^1Z^3=0$.}, i. e.
\be
\label{conic_B}
{\mathcal B}=\{[t^2, 2st, s^2], \quad [s, t]\in \CP^1\}.
\ee
\begin{lemma}
\label{lemma1}
Any point $p\in \CP^3\setminus{\mathcal C}$ lies on a unique secant of ${\mathcal C}$, and thus gives a pair of points on ${\mathcal C}$ (or a tangent in a 
limiting case).
These are the roots of the quadric $Q(p)$.
\end{lemma}
\noindent
{\bf Proof.}
Given $p$ not in ${\mathcal C}$, 
seek two points $(u, v)\in {\mathcal C}$ such that
$(p, u, v)$ lie on the same line in $\CP^3$. Therefore
\be
\label{secant}
p_{ABC} =tu_A u_B u_C+sv_Av_Bv_C
\ee
and $Q(p)_{AB}=tsu_{(A}v_{B)} (u\cdot v)^2$ which is a quadratic
with roots $(u, v)$. The points $(u, v)$ necessarily exist, as given $p$ the
equation (\ref{secant}) is a system of four equations with four unknowns
$(t, s, u, v)$. Tangents and secants of ${\mathcal C}$ are pairwise disconnected, so $p$ belongs to a unique secant.
For the uniqueness consider
\[
u_Au_Bu_C+t(v_Av_Bv_C-u_Au_Bu_C)=w_Aw_Bw_C+s(k_Ak_Bk_C-w_Aw_Bw_C),
\]
and contract both sides with $u^{(A}v^Bw^{C)}$. This gives $(k\cdot w)(k\cdot u) (k\cdot w)s=0$.
So (if all points are distinct) $s=0$, and then $t=0$.
  \begin{center}
  \includegraphics[width=7cm,height=5cm,angle=0]{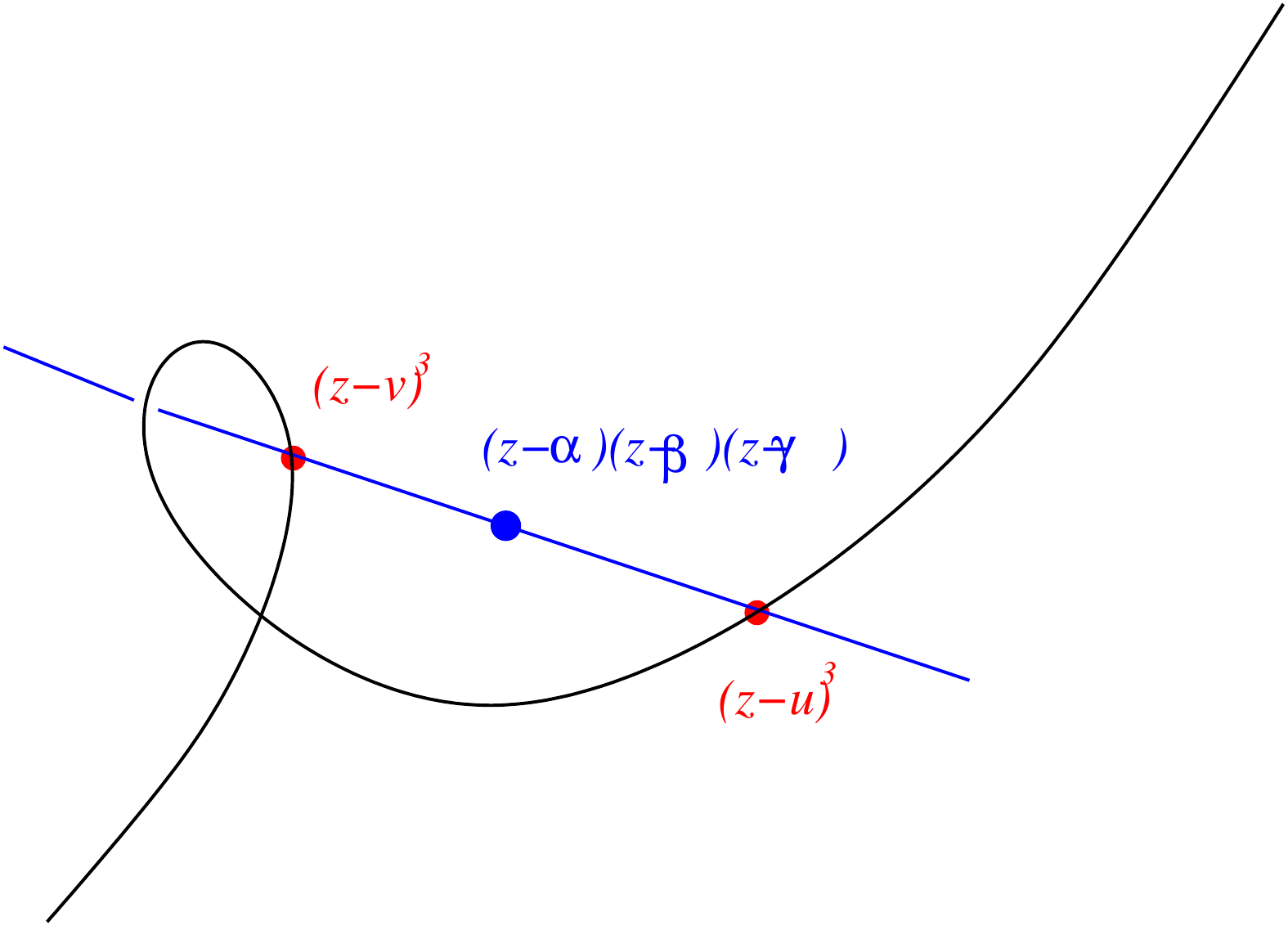}
  \begin{center}
  {{\bf Figure 1.} {\em Secant of the twisted cubic}}
  \end{center}
  \end{center}
 \koniec
In Lemma \ref{lemma2} we shall
regard $V_3$ as the $4$--dimensional symplectic representation space
of $SL(2, \C)$.
\begin{lemma}
\label{lemma2}
The map (\ref{map}) is the projectivisation of the moment map for the symplectic $SL(2, \C)$ action
on the space of cubics $\C^4=\mbox{Sym}^{3}(\C^2)$.
\end{lemma}
\noindent
{\bf Proof.}
Let $p_{ABC}=p_{(ABC)}$ be coordinates on $\C^4$. The $SL(2, \C)$--invariant
symplectic form on $\C^4$ is
\[
{\bf \Omega}=dp_{ABC}\wedge dp^{ABC}.
\]
It is preserved by the action generated by the three vector fields
\be
\label{ham_vector}
H_{AB}=2{p_{(A}}^{CD}\p_{B)CD}, \quad\mbox{where}\quad \p_{BCD}=\frac{\p}{\p p^{BCD}}.
\ee
Therefore
\begin{eqnarray*}
H_{AB}\hook{\bf\Omega}&=&d ({p_{A}}^{DE}p_{BDE})\\
&=&d (\xi_{AB}). 
\end{eqnarray*}
This
extends to the projectivisaton, where $[\xi_{00}, \xi_{01}, \xi_{11}]\in \CP^2$
are homogeneous coordinates of a point $Q(p)$,
where the quadratic map
$Q:\C^4\rightarrow \mathfrak{sl}(2)$ is the projectivisation of the 
moment map.
\koniec
 There are three orbits of the $SL(2)$ symplectic action on $\CP^3$. The generic orbit
corresponding to cubics with three distinct roots, the one dimensonal orbit ${\mathcal C}$, and
the two dimensional orbit ${\mathcal B}$ of cubics with two roots. The union of ${\mathcal C}$
and $\mathcal{B}$ is the discriminant divisor which meets any line in four points.
This divisor is a quartic surface in $\CP^3$ - the union of all tangents 
to $\mathcal{C}$.

\subsection{The Gergonne conic}
Let $p(z)=(z-\alpha)(z-\beta)(z-\gamma)$ be a cubic with three distinct roots\footnote{In this section we shall use $z$ as an affine coordinate on $\CP^1$. Thus if $z^A=[x, y]$ and $\alpha_{A}=[\alpha_0, \alpha_1]$
then $<z, \alpha>_1=(z-\alpha)$, where $z=x/y$ and $\alpha=-\alpha_1/\alpha_0$.}.
Three quadratics
with roots $(\alpha, \beta),  (\alpha, \gamma)$ and $(\beta, \gamma)$ correspond to vertices
of a triangle $T_{(\alpha, \beta, \gamma)}$ in $\CP^2$. This triangle is circumscribed about the conic ${\mathcal B}\subset \CP^2$ given by (\ref{conic_B}).
 
The triangle $T_{(\alpha, \beta, \gamma)}$ is tangent to ${\mathcal B}$ at three points
corresponding to quadratics $(z-\alpha)^2, (z-\beta)^2$ and $(z-\gamma)^2$.
Connect each point of tangency with the opposite vertex of the triangle by a 
line,  i. e. $(z-\alpha)^2$ connects to $(z-\beta)(z-\gamma)$ etc. The limiting
case of the Brianchon's 
theorem\footnote{Brianchon's theorem is a converse to Pascal's theorem in projective geometry. It states that principal diagonals of a hexagon circumscribed around a conic section meet at a single point. In the limiting
case one edge of the hexagon degenerates to a point, and the opposite
three edges degenerate to a segment of a line. In this limit the hexagon becomes a triangle.}  implies that the resulting three
lines intersect at one point. If the conic ${\mathcal B}$ were a circle, $Q(p)$ would be
the Gergonne point of the triangle $T_{(\alpha, \beta, \gamma)}$.
We shall call $Q(p)$ the Gergonne point in general.

The map (\ref{map}) singles out a point in $\CP^2$, and the following Lemma
shows that this is the Gergonne point of the triangle 
$T_{(\alpha, \beta, \gamma)}$.
\begin{lemma}
The Gergonne point is the image of $p$ under the map (\ref{map}).
\end{lemma}
\noindent
{\bf Proof.}
  \begin{center}
  \includegraphics[width=7cm,height=5cm,angle=0]{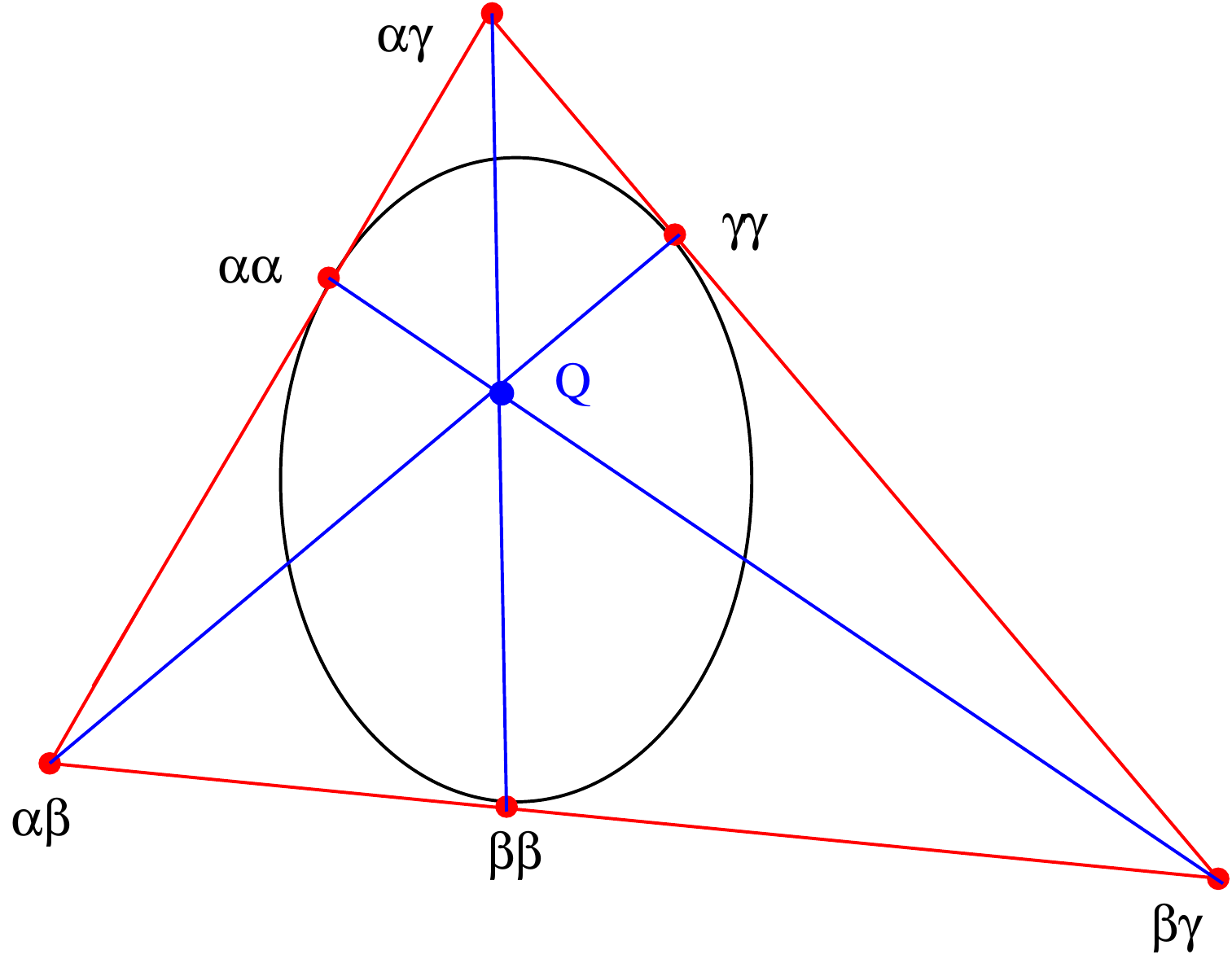}
  \begin{center}
  {{\bf Figure 2.} {\em The Gergonne point}}
  \end{center}
  \end{center}
Consider a line joining the 
vertices of $T_{(\alpha, \beta, \gamma)}$ to the opposite tangency points
\begin{eqnarray*}
L_1(t_\gamma)&=&(z-\alpha)(z-\beta)+t_\gamma ((z-\alpha)(z-\beta)-(z-\gamma)^2)\\
L_2(t_\beta)&=&(z-\alpha)(z-\gamma)+t_\beta ((z-\alpha)(z-\gamma)-(z-\beta)^2)\\
L_3(t_\alpha)&=&(z-\gamma)(z-\beta)+t_\alpha ((z-\gamma)(z-\beta)-(z-\alpha)^2).
\end{eqnarray*}
The system of equations $L_1(t_\gamma)=L_2(t_\beta)=L_3(t_\alpha)$
admits a unique solution for $(t_\alpha, t_\beta, t_\gamma)$ so the Gergonne point indeed exists, and the 
corresponding element of $V_2$ is a quadric with the roots\footnote{These two roots are the fixed points of the Mobius transformation which permutes the roots of the cubic: 
$(\alpha, \beta, \gamma)\rightarrow (\gamma, \alpha, \beta)$.}
\be
\label{greg_guadric}
\frac{\alpha^2(\beta+\gamma)+\beta^2(\alpha+\gamma)+\gamma^2(\alpha+
\beta)-6\alpha\beta\gamma\pm\sqrt{3}i(\alpha^2(\beta-\gamma)+
\beta^2(\gamma-\alpha)+\gamma^2(\alpha-\beta))}{(\alpha-\beta)^2+(\alpha-\gamma)^2+(\beta-\gamma)^2}.
\ee
Now compute $Q(p)$, where $p=(z-\alpha)(z-\beta)(z-\gamma)$, and find
that the roots of the resulting quadric conicide with (\ref{greg_guadric}).
\koniec

\vskip3pt
  Consider a line $L_{(p, q)}\subset \CP^3$ containing two distinct 
cubics $p$ and $q$ not in ${\mathcal C}$. This is
given by
\[
L_{(p, q)}=tp+sq,
\]
where the cubic on the RHS has roots $\alpha(s, t), \beta(s, t), \gamma(s, t)$.
The corresponding conic in $\CP^2$ is 
\be
\label{conic_from_c}
Q(L_{(p, q)})=t^2<p, p>_2+2st<p, q>_2+s^2<q, q>_2.
\ee
This conic intersects the branch conic ${\mathcal B}$ at four points which correspond to 
zeros of the quartic (in $[t, s]$) 
$<Q(L_{(p, q)}), Q(L_{(p, q)})>_2.$ The conic $Q(L)$ is the 
locus of  the Gergonne points of triangles 
$T_{(\alpha(s, t), \beta(s, t), \gamma(s, t))}$. Let us call $Q(L)$ the Gergonne 
conic.
\subsubsection{Characterisation of Gergonne conics}
We want to characterise the conics of the form 
(\ref{conic_from_c}) as a hypersurface in the space $\PPP^5$ of all conics.
Attempting to do it by brute force
leads to  $9$ quadratic equations 
\[
Q(L_{(p, q)})=[\xi_1 t^2+2\xi_2 st+\xi_3 s^2, \xi_4 t^2+2\xi_5 st+\xi_6 s^2,
\xi_7 t^2+2\xi_8 st+\xi_9 s^2]
\]
for
$8$ coefficients $(p^{ABC}, q^{ABC})$, and in principle a sequence of resultants should give a condition.  We shall instead make use of the isomorphism
\be
\label{isomorphism}
\mbox{Sym}^2(\C^3) = \sym^2(\sym^2(\C^2)) = \sym^4(\C^2) \oplus \sym^0(\C^2),
\ee
and express  quadratic forms on $\C^3$ as pairs consisting of a binary quartic, and a scalar. To make (\ref{isomorphism}) explicit
set
\be
\label{the_zets}
Z^1=\frac{1}{\sqrt{2}}(\xi^{00}+\xi^{11}), \quad
Z^2=-\frac{i}{\sqrt{2}}(\xi^{00}-\xi^{11}), \quad Z^3=i\sqrt{2}\xi^{01}
\ee
so that $2Z^TZ=<\xi, \xi>_2$. 
\vskip5pt
Let $\beta\in \sym^2(\sym^2(\C^2))$, 
so that in spinor notation
\[
\beta_{ABCD}=\beta_{(AB)CD}=\beta_{AB(CD)}=\beta_{CDAB}.
\]
Let $\xi^{AB}=\xi^{(AB)}$ be homogeneous coordinates on $\CP^2$. Any conic then
takes the form
\[
\beta_{ABCD}\xi^{AB}\xi^{CD}=0,
\]
or
\begin{eqnarray}
\label{conic_2}
&&24f_0 (\xi^{11})^2-96 f_1 \xi^{01}\xi^{11}+48f_2\xi^{00}\xi^{11}+96f_2(\xi^{01})^2
-96f_3\xi^{00}\xi^{01}\nonumber\\
&&+24f_4 (\xi^{00})^2+4G(\xi^{00}\xi^{11}-(\xi^{01})^2)=0
\end{eqnarray}
where
\[
\F=f_0x^4+4f_1x^3y+6f_2x^2y^2+4f_3xy^3+f_4y^4
\]
is the binary quartic corresponding to $\beta_{(ABCD)}$ and $G$ is a multiple
of ${\beta_{AB}}^{AB}$. The equation (\ref{conic_2}) can be equivalently written
as
\be
\label{conic_3}
<\F,\xi^2>_4+G<\xi, \xi>_2=0,
\ee
where $\xi^2=<\xi, \xi>_0$.
\begin{prop}
\label{propFG}
Let $(\F, G)\in \sym^{4}(\C^2)\oplus \sym^0(\C^2)$ represent a conic in $\CP^2$.
This conic is an image of a line in $\CP^3$ under the map (\ref{map}) iff
\be
\label{invariant_1}
{\mathcal I}:= 4 <\F,\F>_4 -  G^2 = 0.
\ee
\end{prop}
\noindent
{\bf Proof.}
We aim to characterise conics of the form  (\ref{conic_from_c})  as a 
hypersurface in the space of all conics. One can verify (using MAPLE) the following identity\footnote{We are grateful to Robert Bryant for pointing out this identity to us.}
\[
< <p,q>_1 , Q(tp+sq)^2 >_4 
+ 6 <p, q>_3 < Q(tp+sq), Q(tp+sq) >_2 = 0,
\]
where $<p ,q>_1$ is nonzero as long as $p$ and 
$q$ are linearly independent in $V_3$.
Therefore the conic in $\CP^2$ traced out by
$Q(tp+sq)$ 
satisfies the quadratic equation
\[
< <p, q>_1 , <Q(tp+sq),Q(tp+sq)>_0 >_4   + 6 <p ,q>_3 <Q(tp+sq), Q(tp+sq)>_2 = 0.
\]
Thus, it follows that the equation of the Gergonne conics 
is given by pairs 
\[( <p, q>_1, 6 <p, q>_3 ) \in \sym^4(\C^2) \oplus \sym^0(\C^2).
\]
The pairs $(\F,G) = ( <p ,q>_1, 6 <p ,q>_3 )$
satisfy the quadratic relation (\ref{invariant_1})
so that this equation defines the hypersurface in the space of conics.
\koniec
We shall now replace (\ref{invariant_1}) by a condition on the $3$ by $3$ symmetric matrix
representing a conic.
\begin{prop}
\label{prop_matrix}
Let $A=A^T\in \sym^2(\C^3)$ represent a conic 
\be \label{conic_4}
Z A Z^T=0,
\ee where
$Z=[Z^1, Z^2, Z^3]\in \CP^2$. This conic is the image of a line in $\CP^3$ under
(\ref{map}) iff
\be
\label{inv_2}
{\mathcal I}:=\frac{1}{6}\mbox{Tr}(A^2)-\frac{1}{12}\mbox{Tr}(A)^2=0.
\ee
\end{prop}
\noindent
{\bf Proof.}
Consider (\ref{conic_4}), and substitute the expressions (\ref{the_zets}) for
$Z^{i}$. Comparing (\ref{conic_4}) and (\ref{conic_2}) yields
\begin{eqnarray*}
f_0&=&\frac{1}{48}(A_{11}-A_{22}+2iA_{12}), \quad 
f_1 =
\frac{1}{48}(A_{23}-iA_{13}),\quad
f_2=\frac{1}{144}(A_{11}+A_{22}-2iA_{33}), \\
f_3 &=&
-\frac{1}{48}(A_{23}+iA_{13}), \quad
f_4=\frac{1}{48}(A_{11}-A_{22}-2iA_{12}), \quad G=\frac{1}{6}\mbox{Tr}(A).
\end{eqnarray*}
We now use the result of Propositon \ref{propFG}, and verify
that
\begin{eqnarray*}
{\mathcal I}&=&4<\F, \F>_4-G^2\\
&=&\frac{1}{12}(A_{11}^2+A_{22}^2+A_{33}^3+4A_{12}^2+4A_{13}^2+
4A_{23}^2-2A_{11}A_{22}-2A_{11}A_{33}-2A_{22}A_{33}\\
&=&\frac{1}{6}\mbox{Tr}(A^2)-\frac{1}{12}\mbox{Tr}(A)^2.
\end{eqnarray*}
\koniec
We are now ready to establish Theorem \ref{theo3_intro} from the Introduction
\begin{theo}
\label{theo_ponce}
Let the conic ${\mathcal A}$  given by (\ref{conic_3}) 
be the image of a line in $\CP^3$ under the map
(\ref{map}). Then the matrices $(A, {\bf I})$  form a Poncelet pair 
with a triangle\footnote{Hitchin argues that the vertices of all triangles 
$T_{(\alpha(s, t), \beta(s, t), \gamma(s, t))}$ sweep
a conic - call this one the Hitchin conic $H(L)$ - so that the Hitchin conic, and the branch conic ${\mathcal B}$ form a Poncelet pair. On page $17$ in \cite{poncelet_nigel} he takes
\[
tp(\alpha(s, t))+sq(\alpha(s, t))=0, \quad tp(\beta(s, t))+sq(\beta(s, t))=0
\]
so that
\[
0=p(\alpha(s, t)) q(\beta(s, t))  -q(\alpha(s, t)) p(\alpha(s, t))
=<\alpha(s, t), \beta(s, t)> R(\alpha(s, t), \beta(s, t))
\]
where $R\in(\sym^4(\C^2)\oplus\sym^0(\C^2)$ is  polynomial in $\alpha, \beta$, which defines a conic $H(L)$. 
We claim that the Hitchin conic does not coincide with the Gergonne conic.
To see this consider the Gergonne conic (\ref{conic_from_c})
corresponding to cubics $p=(z-a)(z-b)(z-c)$ and $q=(z-A)(z-B)(z-C)$.
Now consider the vertices of the triangles $T_{(abc)}$ and $T_{(ABC)}$, i. e.
quadrics $(z-a)(z-b)$ etc. Do they belong to the Gergonne conics? If they do, then
the Hitchin and Gergonne conics must coincide, as five points determine a conic.
It is sufficient to verify it for one vertex, say $(z-a)(z-b)$ - the others will follow by symmetry. Consider the transvectant $<Q(L), (z-a)^2>_2$. This gives
a quadratic in $[t, s]$ with two roots. Do the same with $(z-b)^2$ - another two roots. Computation shows that there is no pair of common roots,
so $(z-a)(z-b)$ does not belong to $Q(L)$ for any $[t, s]$.
A combination of Hitchin's porism with our result shows that
the locus
of Gergonne points of a poristic family of triangles is itself 
a conic. If the conics are coaxial circles, this was noted
in \cite{Gallaty}.}: there exists a triangle  inscribed in 
${\mathcal A}$ and  circumscribed about the
base conic
${(Z^1)}^2+{(Z^2)}^2+{(Z^3)}^2=0$.
\end{theo}
\noindent
{\bf Proof.}
Let $A$ and $D$ be symmetric matrices defining two conics 
${\mathcal A}$ and ${\mathcal D}$. Cayley \cite{cayley,gh} gave an algebraic
conditions for $(A, D)$ in order for a polygon with $N$ vertices
 inscribed in ${\mathcal A}$ and  circumscribed about ${\mathcal D}$  exist.
Consider the expansion
\[
\sqrt{\mbox{det}(sA-D)}=a_0+s a_1+s^2a_2+\dots
\]
Assume that $N=2n+1$ is odd (a similar formula exists for even $N$). Then the 
necessary and sufficient condition for the existence of an $N$--gon is
\[
\begin{vmatrix}
a_2 & \dots & a_{n+1} \\ 
. & .     & . \\ 
. & .  & . \\ 
a_{n+1} & \dots & a_{2n}  \notag
\end{vmatrix}
=0.
\]
We need to consider the special case when $N=3$, and $D=\mbox{diag}(1, 1, 1)$.
The characteristic polynomial for $3$ by $3$ matrices with unit determinant is
\[
\mbox{det}(sA-{\bf I})=1+s\mbox{Tr}(A)+\frac{1}{2}s^2(\mbox{Tr}(A)^2-\mbox{Tr}(A^2))
+s^3.
\]
Comparing this with
\[
\sqrt{\mbox{det}(sA-{\bf I})}=1+sa_1+s^2a_2+\dots
\]
gives 
\[
a_1=\frac{1}{2}\mbox{Tr}(A), \quad a_2=\frac{1}{8}(\mbox{Tr}(A)^2-2\mbox{Tr}(A^2)).
\]
The vanishing of $a_2$ is the Cayley
condition\footnote{Robert Bryant has suggested that if the projection (\ref{map}) is replaced by a generic triple of quadric, 
the relation (\ref{invariant_1}) becomes quartic (and not quadratic). If this quartic is $SO(3)$ invariant, then it must take the form
\[
c_0\mbox{Tr}(A^3)\mbox{Tr}(A)+c_1\mbox{Tr}(A^2)(\mbox{Tr}(A))^2+
c_2(\mbox{Tr}(A^2))^2+c_3(\mbox{Tr}(A))^4.
\]
We have checked that (dehomogenising with $\mbox{det}(A)$) this quartic does not satisfy the system
(\ref{moraru_eq}) for any choice of the $c_i$. Thus if there is a quartic relation, then it is not expressible
by traces.} for $(A, {\bf I})$ to form a Poncelet pair with a triangle, and
we see that $a_2$ is a constant multiple of (\ref{inv_2}). The result now follows from Proposition \ref{prop_matrix}.
\koniec
\section*{Appendix}
\appendix
\setcounter{equation}{0}
\def\theequation{\thesection{A}\arabic{equation}}
\subsection*{Two--component spinors}
\label{sec_spinors}
A convenient way to represent binary quartics and the associated
invariants uses the two--component spinor notation \cite{md_rp}.

Let the capital letters $A, B, \dots$ denote indices taking values
$0$ and $1$. The general quartic is represented by a 
{\it symmetric spinor} of valence $4$. 
Let $\pi_{A}=[s, t]$ be homogeneous coordinates on  $(\C^2)^*$
(the dual of $\C^2$).
A homogeneous 
quartic corresponding to a vector $V\in TM$ is given by 
\[
V=V_{ABCD}\pi^A\pi^B\pi^C\pi^D=\alpha t^4+4\beta t^3s+6\gamma
t^2s^2+4\delta ts^3+\epsilon s^4,
\] 
where $\pi^A=[t, -s]$. The spinor indices are lowered by the anti--symmetric matrix 
$\varepsilon_{AB}$ with $\varepsilon_{01}=1$, so that
\[
\psi_{AB\dots C}=\psi^{PQ\dots R}\varepsilon_{PA}\varepsilon_{QB}\dots
\varepsilon_{RC}.
\]
Let $e^1, \dots, e^5$ be one--forms on $M$. Define $e^{ABCD}$ by 
\be
\label{abcd}
e^1=e^{1111}, \quad e^{2}=e^{0111}, \quad e^{3}=e^{0011},\quad
e^4=e^{0001}, \quad e^5=e^{0000}.
\ee
The $GL(2)$--structure on $M$ is given by $S=\pi_{A}\pi_B\pi_C\pi_De^{ABCD}$,
and the $SO(3)$ structure (\ref{so3vector}) is given by
\[
g=e^{ABCD}\odot e_{ABCD}, \quad G={e^{AB}}_{CD}\odot {e^{CD}}_{EF}\odot {e^{EF}}_{AB}.
\]
\subsection*{Miscellanous formulae}
Let $(a, b, p, q, r)$ be local coordinates on $M=SL(3, \R)/SO(3)$ such that
a symmetric matrix $A\in M$ with $\mbox{det}(A)=1$ is given by
\[
A=BB^T, \quad \mbox{where}\quad
B=\left(
\begin{array}{ccc}
e^{-a-b} & pe^b& qe^a\\
0 & e^b &re^a\\ 
0 & 0 &e^a
\end{array}
\right ).
\]
The Riemannian Einstein metric on $M$ is given by
\begin{eqnarray*}
g&=&4\; \mbox{Tr}(A^{-1}dA\cdot A^{-1}dA)\\
&=&8[(2da+db)^2+ 3db^2+
e^{2a+4b}dp^2+e^{2a-2b}dr^2+e^{4a+2b}(dq-pdr)^2].
\end{eqnarray*}
The $8$ dimensional isometry group $SL(3, \R)$ of $g$ is generated by the Killing vector fields
\begin{eqnarray}
X_1&=&\p_p+r\p_q, \quad X_2=\p_q, \quad X_3=\p_r,\nonumber\\
X_4&=&\p_a-p\p_p-2q\p_q-r\p_r, \quad X_5=\p_b-2p\p_p-q\p_q+r\p_r,\nonumber\\
X_6&=&p\p_b-(1+p^2-e^{-2a-4b})\p_p-r\p_q+q\p_r,\nonumber \\
X_7&=&r\p_a-r\p_b  -(1+r^2-e^{2b-2a})\p_r+(e^{2b-2a}p-rq)\p_q+(pr-q)\p_p\nonumber\\
X_8&=&q\p_a-rp\p_b+(p^2r-re^{-2a-4b}-qp)\p_p+
(p^2e^{2b-2a}+e^{-4a-2b}-q^2-1)\p_q+(e^{2b-2a}p-rq)\p_r.\nonumber
\label{5d_killing}
\end{eqnarray}
The group of rotations $SO(3)\subset SL(3, \R)$ is generated by
${X_6, X_7, X_8}$ as
\[
[X_6, X_7]=X_8, \quad [X_6, X_8]=-X_7, \quad [X_7, X_8]=X_6.
\]
There are two invariants of the $SO(3)$ action on $M$ which we chose to be 
$\mbox{Tr}(A)$ and $\mbox{Tr}(A^2)$. We have verified that there is no non--zero function of these which satisfies the system
(\ref{moraru_eq}).

Let $(M, g, G)$ be the $SO(3)$ structure from Proposition \ref{so3prop}.
The associated differential operators which characterise the range of the Penrose--Radon transform in Theorem \ref{theo_radon} are given by
\be
\label{1st_operator}
F\rightarrow (\Delta_g+\frac{1}{12}\mbox{Id}\Big)F, \quad\mbox{where}
\ee
\begin{eqnarray}
\label{aplaplacian}
\Delta_g&=&g^{ab}\nabla_a\nabla_b\nonumber\\
&=&\frac{1}{24}\Big(
\frac{\p^2}{\p a^2}+\frac{\p^2}{\p b^2}
-\frac{\p^2}{\p a\p b}+3\frac{\p}{\p a}\Big)\\
&+&\frac{1}{8}e^{2b-2a}\Big(p^2\frac{\p^2}{\p q^2}+
\frac{\p^2}{\p r^2}+2p\frac{\p^2}{\p r\p q} 
 +e^{-2a-4b}\frac{\p^2}{\p q^2}
+e^{-6b}\frac{\p^2}{\p p^2}
\Big)\nonumber,
\end{eqnarray} 
and
\be
\label{2nd_operator}
F\rightarrow {G^{ab}}_c\nabla_a\nabla_b F-\frac{1}{24}\nabla_c F.
\ee
A function $F:M\rightarrow \R$ is in the kernel of the operator 
(\ref{2nd_operator}) iff
\begin{eqnarray}
\label{big_operator}
&&
\Big(-3e^{2b-2a}p^2\frac{\p^2}{\p q^2}-6e^{2b-2a}p\frac{\p^2}{\p r\p q}
+3e^{-2a-4b}\frac{\p^2}{\p p^2}-3e^{2b-2a}\frac{\p^2}{\p r^2}
+2\frac{\p^2}{\p a \p b}-\frac{\p^2}{\p a^2}+3\frac{\p}{\p b}-2\frac{\p}{\p a}
\Big)Fda\nonumber\\
&&
+
\Big(-3e^{2b-2a}p^2\frac{\p^2}{\p q^2}-6e^{2b-2a}p\frac{\p^2}{\p r\p q}
+3e^{-4a-2b}\frac{\p^2}{\p q^2}-3e^{2b-2a}\frac{\p^2}{\p r^2}-\frac{\p^2}{\p b^2}
+2\frac{\p^2}{\p a\p b}+\frac{\p}{\p b}
\Big)Fdb\nonumber\\
&&
-\Big(  
3e^{2b-2a}p\frac{\p^2}{\p q^2}+3e^{2b-2a}\frac{\p^2}{\p r\p q}
+\frac{\p^2}{\p b\p p}-2\frac{\p^2}{\p p\p a}-\frac{\p}{\p p}
\Big)Fdp\nonumber\\
&&
-\Big(3p\frac{\p^2}{\p q\p p}+\frac{\p^2}{\p q\p a}-2\frac{\p^2}{\p q \p b}+3\frac{\p^2}{\p r\p p}+2\frac{\p}{\p q}
\Big)Fdq\nonumber\\
&&
-\Big(-3(p^2-e^{-2a-4b})\frac{\p^2}{\p q\p p}+
3p\frac{\p^2}{\p q\p b}-3p\frac{\p^2}{\p r\p p}+\frac{\p^2}{\p r\p b}
+\frac{\p^2}{\p r \p a}+2\frac{\p}{\p r}
\Big)Fdr=0.\nonumber
\end{eqnarray}
\subsection*{The generalised Legendre transform.}
In this paper we have studied the interplay between
the $SO(3)$ structure on $M=SL(3, \R)/SO(3, \R)$, the Penrose--Radon transform on conics, 
and the tri--K\"ahler metrics whose twistor spaces admit a  holomorphic fibration
over $\CP^2$. In this Appendix we shall put an analogous,  but  simpler construction
of the generalised Legendre transform \cite{IR96,B00,DM} in this framework. Here $M=\R^5$, the
$SO(3)$--structure consists of a flat metric and a three--form preseved by an $8$--dimensional subgroup of the full isometry group, and the twistor
space is an affine bundle over the total space of $\OO(4)\rightarrow \CP^1$.
Let ${\mathcal Z}\rightarrow \OO(4)$ by one such affine bundle corresponding to a cohomology class $[f]\in H^1(\OO(4), \OO_{\CP^1}(-2))$. The sections of $\OO(4)\rightarrow\CP^1$ are binary quartics parametrised by points in 
$M_\C=\mbox{Sym}^{4}(\C^2)=\C^5$. Let $L_m$ be a section of $\OO(4)$ corresponding
to a point $m\in M_\C$.
The {\em real} sections are preserved by an anti--holomorphic involution on $\OO(4)$ which restricts to an antipodal map on each section. They are of the form
\be
\label{omega_coord}
[s, t]\rightarrow \omega=s^4\tau+ts^3 z+t^2s^2 y-t^3s\ov{z}+t^4\ov{\tau},
\ee
where $[s, t]$ are homogeneous coordinates on $\CP^1$, and $(z, \ov{z}, \tau, 
\ov{\tau}, y)$ are coordinates on $M=\R^5$. The $GL(2)$ structure (\ref{para_con_exp}) on $M_\C$
is given by ${T}_mM_\C\cong  H^0(L_m, \OO(4))=\mbox{Sym}^4(\C^2)$, which restricts
to 
\[
S=  s^4d\tau+ts^3 dz+t^2s^2 dy-t^3sd\ov{z}+t^4d\ov{\tau},
\]
on $M$. This gives rise to an $SO(3)$ structure (\ref{so3structure}), where the metric $g$ and the three--form
$G$ are given by (\ref{so3vector}) with
\[
e^1=d\ov{\tau}, \quad e^2=-\frac{1}{4}d\ov{z},\quad e^3=\frac{1}{6}dy, \quad e^4=\frac{1}{4}
dz, \quad e^5=d\tau.
\]
The range of the Penrose--Radon transform
\[
F(z, \ov{z}, \tau, \ov{\tau}, y)=\oint_{\Gamma\subset \CP^1} f(\omega, [s, t]) (sdt-tds)
\]
(where $\omega$ is given by (\ref{omega_coord}))
is characterised by the over--determined system of linear PDEs
\be
\label{legeqq}
\Delta_g F=0, \quad G^{ab}_c\nabla_a\nabla_b F=0
\ee
(note that $\kappa$ and $R$ in (\ref{G5}) and (\ref{G6}) both vanish in this case), or explicitly 
\begin{eqnarray*}
&&F_{y\tau}-F_{zz}=0, \quad F_{\tau\ov{z}}+F_{yz}=0, \quad
F_{\tau\ov{\tau}}+F_{z\ov{z}}=0,\\
&& F_{z\ov{z}}+F_{yy}=0,\quad F_{\ov{\tau}z}+F_{\ov{z}y}=0, \quad F_{\ov{\tau}y}-F_{\ov{z}\ov{z}}=0.
\end{eqnarray*}
Assuming that the $J$--invariant (\ref{j_invariant}) of the binary quartic corresponding to
$dF$ is not zero,   we can now apply Theorem \ref{theo_2_forms} to construct the tri--K\"ahler structure on the four--fold $X=\{m\in M, F(m)=0\}$. 
The metric $g$ on $M$ is flat, and ${\Gamma^A}_B=0$. The identity (\ref{identity}) implies that the two--forms
$\Sigma^{AB}$ are closed iff equations (\ref{legeqq}) hold.
The pull--backs
of $\Sigma^{AB}$ to $M$ are also normalised such that $\Sigma^{(AB}\wedge\Sigma^{BC)}=0$
and of constant length. Thus they define a hyper--K\"ahler structure. 

The original 
set-up of \cite{IR96} would lead to a K\"ahler potential as follows. Let $H:M\rightarrow \R$
be a function such that $F=\p H/\p y$. Now perform the Legendre transform
$u=\p H/\p z$, and eliminate the cordinates $(z, \ov{z}, y)$, using $(\tau, u,  \ov{\tau},\ov{u})$ as holomorphic and antiholomorphic coordinates on $X$.
The Kahler potential
$
{\mathcal K}(\tau, u,  \ov{\tau},\ov{u})= H-\tau u-\ov{\tau}\ov{u}
$
satisfies 
\[ 
{\mathcal K}_{\tau\ov{\tau}}{\mathcal K}_{u\ov{u}}- 
 {\mathcal K}_{\tau\ov{u}}{\mathcal K}_{u\ov{\tau}}=1,
\]
and the barycentre metric 
\[
\gamma_B={\mathcal K}_{\tau\ov{\tau}}d\tau d\ov{\tau}+
{\mathcal K}_{u\ov{u}}du d\ov{u}+
{\mathcal K}_{\tau\ov{u}}d\tau d\ov{u}+
{\mathcal K}_{u\ov{\tau}}du d\ov{\tau}
\]
on $X$ is hyper-K\"ahler.


\end{document}